\documentclass[11pt]{article}
\usepackage{amsmath}
\usepackage{amssymb,amsmath}
\usepackage{mathrsfs}
\usepackage{graphicx}
\DeclareMathOperator{\intd}{d}
\let\pa=\partial
\let\al=\alpha

\let\f=\frac
\let\p=\psi
\let\om=\omega

\def\cC{{\cal C}}

\def\cF{{\cal F}}

\def\cH{{\cal H}}

\def\cP{{\cal P}}

\def\cS{{\cal S}}

\def\cZ{{\cal Z}}

\def\no{\noindent}
\def\na{\nabla}

\def\p{\partial}

\def\dv{\mbox{div}}

\def\eqdefa{\buildrel\hbox{\footnotesize def}\over =}

\def\C{\mathop{\bf C\kern 0pt}\nolimits}
\def\DD{\mathop{\bf D\kern 0pt}\nolimits}
\def\K{\mathop{\bf K\kern 0pt}\nolimits}
\def\N{\mathop{\bf N\kern 0pt}\nolimits}
\def\Q{\mathop{\bf Q\kern 0pt}\nolimits}
\def\R{\mathop{\bf R\kern 0pt}\nolimits}

\newcommand{\Z}{{\mathbf Z}}

\newcommand{\ef}{ \hfill $ \blacksquare $ \vskip 3mm}

\newcommand{\beq}{\begin{equation}}
\newcommand{\eeq}{\end{equation}}
\newcommand{\ben}{\begin{eqnarray}}
\newcommand{\een}{\end{eqnarray}}
\newcommand{\beno}{\begin{eqnarray*}}
\newcommand{\eeno}{\end{eqnarray*}}

\renewcommand{\theequation}{\thesection.\arabic{equation}}

\newtheorem{Theorem}{Theorem}[section]
\newtheorem{Definition}[Theorem]{Definition}
\newtheorem{Proposition}[Theorem]{Proposition}
\newtheorem{Lemma}[Theorem]{Lemma}

\newtheorem{Remark}[Theorem]{Remark}

\begin{document}
\title{Well-posedness in critical spaces for the compressible Navier-Stokes equations with
density dependent viscosities}
\author{Qionglei Chen $^\dag$  Changxing Miao $^{\dag}$ and Zhifei Zhang
$^{\ddag}$\\[2mm]
{\small $ ^\dag$ Institute of Applied Physics and Computational Mathematics, Beijing 100088, China}\\
{\small E-mail: chen\_qionglei@iapcm.ac.cn  and   miao\_changxing@iapcm.ac.cn}\\[2mm]
{\small $ ^\ddag$ School of  Mathematical Science, Peking University, Beijing 100871, China}\\
{\small E-mail: zfzhang@math.pku.edu.cn}}

\date{25 November, 2008}

\maketitle

\begin{abstract}
In this paper, we prove the local well-posedness in critical Besov spaces for the compressible Navier-Stokes equations with
density dependent viscosities under the assumption that the initial density is bounded away from zero.
\end{abstract}

 \begin{minipage}{120mm}

 \end{minipage}
\renewcommand{\theequation}{\thesection.\arabic{equation}}
\setcounter{equation}{0}
%%%%%%%%%%%%%%%%%%%%%%%%%%%%%%%%%%%%%%%%%%%%%%
%%%%%%%%%%%%%%%%%%%%%%%%%%%%%%%%%%%%%%%%%%
\section{Introduction}

In this paper, we consider the compressible Navier-Stokes equations
with density dependent viscosities in $\R^+\times \R^N(N\ge 2)$:
\begin{equation}\label{eq:cNS}
\left\{
\begin{array}{ll}
\p_t\rho+\textrm{div}(\rho u)=0,\\
\p_t(\rho u)+\textrm{div}(\rho u\otimes u)-\textrm{div}(2\mu(\rho)D(u))-\na(\lambda(\rho) \textrm{div}u)+\na P(\rho)=0, \\
(\rho,u)|_{t=0}=(\rho_0,u_0).
\end{array}
\right.
\end{equation}
Here $\rho(t,x)$ and $u(t,x)$ are the density and velocity of the
fluid. The pressure $P$ is a smooth function of $\rho$, $D(u)=\f12(\na u+\na u^t)$ is the strain tensor,
the Lam\'{e} coefficients $\mu$ and $\lambda$ depend smoothly on $\rho$ and
satisfy
\ben\label{assu:coeff}
\mu>0\quad \textrm{and} \quad \lambda+2\mu>0,
\een
which ensures that the operator $-\textrm{div}(2\mu(\rho)D\cdot)-\na(\lambda(\rho) \textrm{div}\cdot)$
is elliptic. An important example is included in the system (\ref{eq:cNS}): the
viscous shallow water equations($N=2,\mu(\rho)=\rho,\lambda(\rho)=0$ and $P(\rho)=\rho^2$).

The local existence and uniqueness of smooth solutions for the
system (\ref{eq:cNS}) were proved by Nash \cite{Nash} for smooth
initial data without vacuum. Later on, Matsumura and Nishida\cite{Mat-Nis}
proved the global well-posedness for smooth data close to
equilibrium, see also \cite{Kaz-She} for one dimension. Concerning
the global existence of weak solutions for the large initial data,
we refer to \cite{Ber-Des, BDL,Lions,Mel-Vas}. We may refer to
\cite{BDM, CMZ-SIAM, Wang-Xu} and references therein for the viscous shallow
water equations.

This paper is devoted to the study of the well-posedness of the
system (\ref{eq:cNS}) in the critical spaces. Recently, Danchin has
obtained several important well-posedness results in the critical
spaces for the compressible Navier-Stokes equations
\cite{Dan-inve,Dan-cpde,Dan-NDEA}. To explain the precise meaning of
critical spaces, let us consider the incompressible Navier-Stokes
equations \beno (NS)\quad\left\{
\begin{array}{ll}
\p_tu-\Delta u+u\cdot\na u+\na p=0, \\
\dv u=0.
\end{array}
\right.
\eeno
It is easy to find that if $(u,p)$ is a solution of (NS), then
\ben\label{scaling}
u_\lambda(t,x)\eqdefa \lambda u(\lambda^2 t,\lambda x),\quad
p_\lambda(t,x)\eqdefa \lambda^2 p(\lambda^2 t,\lambda x)
\een
is also a solution of (NS). For the (NS) equations, a functional space $X$ is critical
if the corresponding norm is invariant under the scaling of (\ref{scaling}).
Obviously, $\dot H^{\frac N 2-1}$ is a critical space. Fujita and Kato\cite{Fuj-Kat} proved the well-posedness
of (NS) in $\dot H^{\frac N 2-1}$, see also \cite{Can1, Can2, Mey} and references therein
for the well-posedness in the other critical spaces. For the compressible Navier-Stokes equations,
let us introduce the following transformation
\beno
\rho_\lambda(t,x)\eqdefa  \rho(\lambda^2 t,\lambda x),\quad
u_\lambda(t,x)\eqdefa \lambda u(\lambda^2 t,\lambda x).
\eeno
Then if $(\rho,u)$ solves (\ref{eq:cNS}), so does $(\rho_\lambda, u_\lambda)$ provided
the viscosity coefficients are constants and the pressure law has been changed into $\lambda^2 P$.
This motivates the following definition:
\begin{Definition}
We will say that a functional space is critical with respect to the scaling of the equations
if the associated norm is invariant under the transformation:
\beno
(\rho,u)\longrightarrow (\rho_\lambda,u_\lambda)
\eeno
(up to a constant independent of $\lambda$).
\end{Definition}

A natural candidate is the homogenous Sobolev space $\dot
H^{N/2}\times \bigl(\dot H^{N/2-1}\bigr)^N$, but since $\dot
H^{N/2}$ is not included in $L^\infty$, we can not obtain a
$L^\infty$ control of the density when $\rho_0\in \dot H^{N/2}$.
Instead, we choose the initial data $(\rho_0,u_0)$ for some
$\bar{\rho_0}$ in a critical homogenous Besov spaces: \beno
(\rho_0-\bar{\rho}_0,u_0)\in \dot B^{\f {N} p}_{p,1}\times
\bigl(\dot B^{\f N p-1}_{p,1}\bigr)^N, \eeno since $\dot B^{\f {N}
p}_{p,1}$ is continuously embedded in $L^\infty$.

However, working in the critical spaces, if we deal with the
elliptic operators of the momentum equations as a constant
coefficient second order operator plus a perturbation induced by the
density and viscosity coefficients, the perturbation will be a
trouble term. In the case when $\rho-\bar{\rho}_0$ is small in $\dot
B^{ \frac{N}{p}}_{p,1}$ or has more regularity, the perturbation can
be treated as a harmless source term and the corresponding
local-well posedness can be obtained by following the argument of
Danchin \cite{Dan-cpde}, see \cite{Has}.

The purpose of the present paper is to obtain a local well-posedness
result in the critical Besov spaces under the natural physical
assumption that the initial density is bounded away from zero. Our
new observation is that if $\rho-\bar{\rho}_0$ is small in the
weighted Besov spaces $\dot B^{ \frac{N}{p}}_{p,1}(\om)$(see Section
3 for the definition), the perturbation can still be treated as a
harmless source term. Similar idea has been used by the authors of
this paper to prove the local well-posedness in $\dot
B^{1}_{2,1}\times \bigl(\dot B^{0}_{2,1}\bigr)^2$ for the viscous
shallow water equations \cite{CMZ-SIAM}. Very rencently,
Danchin\cite{Dan-cpde07} proved a similar result for the system
(\ref{eq:cNS}) with constant coefficients. The key of his proof is a
new and interesting estimate for a class of parabolic systems with
the coefficients in $C([0,T]; \dot B^{N/2}_{2,1})$. It seems to be
possible to adapt his method to the present model. Here we would
like to present a general functional framework to deal with the
local well-posedness in the critical spaces for the compressible
fluids.\vspace{0.1cm}

Our main result is as follows:

\begin{Theorem}\label{Them1.2} Let $\bar{\rho}_0$ and $c_0$ be two positive constants.
Assume that the initial data satisfies
\beno
(\rho_0-\bar{\rho}_0,u_0)\in \dot B^{\f {N} p}_{p,1}\times
\bigl(\dot B^{\f N p-1}_{p,1}\bigr)^N\quad \textrm{and} \quad \rho_0\ge c_0.
\eeno
Then there exists a positive time $T$ such that\vskip 0.1cm

(a)\, {\bf Existence:}\, If $p\in (1,N]$, the system (\ref{eq:cNS}) has a solution $(\rho-\bar{\rho}_0,u)\in E^p_T$ with
\beno
E^p_T\eqdefa C([0,T]; \dot B^{\f {N} p}_{p,1})\times \Bigl(C([0,T]; \dot B^{\f {N}
p-1}_{p,1}) \cap L^1(0,T; \dot B^{\f {N} p+1}_{p,1})\Bigr)^N,\quad \rho\ge
\f12 c_0;
\eeno

(b)\, {\bf Uniqueness:} If $p\in (1,N]$, then the uniqueness holds in $E^p_T$.
\end{Theorem}
\begin{Remark}
If the Lam\'{e} coefficients $\mu$ and $\lambda$ are constants
satisfying (\ref{assu:coeff}), then the range of $p$ in the existence
result of the system \eqref{eq:cNS} can be extended to $p\in(1,2N)$,
since we can take $p\in (1, 2N)$ in Proposition \ref{Prop:momenequ} for the case when $\overline{\lambda}$
and $\overline{\mu}$ are constants.
\end{Remark}

The structure of this paper is as follows:\vskip 0.1cm

In Section 2, we recall some basic facts about the Littlewood-Paley
decomposition and the functional spaces. In Section 3, we firstly
introduce the weighted Besov spaces, then present some nonlinear
estimates. Section 4 is devoted to the estimates in the weighted
Besov spaces for the linear transport equation. Section 5 is devoted
to the estimates in the weighted Besov spaces for the linearized
momentum equation. In Section 6, we prove the existence of the
solution. In Section 7, we prove the uniqueness of the solution.

\setcounter{equation}{0}
\section{Littlewood-Paley theory and the functional spaces}

Let us introduce the Littlewood-Paley
decomposition. Choose a
radial function  $\varphi \in {\cS}(\R^N)$ supported in
${\cC}=\{\xi\in\R^N,\, \frac{3}{4}\le|\xi|\le\frac{8}{3}\}$ such
that \beno \sum_{j\in\Z}\varphi(2^{-j}\xi)=1 \quad \textrm{for
all}\,\,\xi\neq 0. \eeno The frequency localization operator
$\Delta_j$ and $S_j$ are defined by
\begin{align}
\Delta_jf=\varphi(2^{-j}D)f,\quad S_jf=\sum_{k\le
j-1}\Delta_kf\quad\mbox{for}\quad j\in \Z. \nonumber
\end{align}
With our choice of $\varphi$, one can easily verify that
\beq\label{orth}
\begin{aligned}
&\Delta_j\Delta_kf=0\quad \textrm{if}\quad|j-k|\ge 2\quad
\textrm{and}
\quad \\
&\Delta_j(S_{k-1}f\Delta_k f)=0\quad \textrm{if}\quad|j-k|\ge 5.
\end{aligned}
\eeq We denote the space ${\cZ'}(\R^N)$ by the dual space of
${\cZ}(\R^N)=\{f\in {\cS}(\R^N);\,D^\alpha \hat{f}(0)=0;
\forall\alpha\in\N^d \,\mbox{multi-index}\}$, it also can be
identified by the quotient space of ${\cS'}(\R^N)/{\cP}$ with the
polynomials space ${\cP}$. The formal equality
\beno
f=\sum_{k\in\Z}\Delta_kf \eeno
holds true for $f\in {\cZ'}(\R^N)$ and is called the homogeneous Littlewood-Paley
decomposition.

The operators  $\Delta_j$ help us recall the definition of the Besov
space(see also \cite{Tri}).

\begin{Definition} Let $s\in\R$, $1\le p,
r\le+\infty$. The homogeneous Besov space $\dot{B}^{s}_{p,r}$ is
defined by
$$\dot{B}^{s}_{p,r}=\{f\in {\cZ'}(\R^N):\,\|f\|_{\dot{B}^{s}_{p,r}}<+\infty\},$$
where \beno
\|f\|_{\dot{B}^{s}_{p,r}}\eqdefa \Bigl\|2^{ks}
\|\Delta_kf(t)\|_{p}\Bigr\|_{\ell^r}.\eeno
\end{Definition}

We next introduce the Besov-Chemin-Lerner space
$\widetilde{L}^q_T(\dot{B}^{s}_{p,r})$ which is initiated in
\cite{Che-Ler}.

\begin{Definition}Let $s\in\R$, $1\le
p,q,r\le+\infty$, $0<T\le+\infty$. The space
$\widetilde{L}^q_T(\dot{B}^s_{p,r})$ is defined as the set of all
the distributions $f$ satisfying
$$\|f\|_{\widetilde{L}^q_T(\dot{B}^{s}_{p,r})}<+\infty,$$\
where
$$\|f\|^r_{\widetilde{L}^q_T(\dot{B}^{s}_{p,r})}\eqdefa \Bigl\|2^{ks}
\|\Delta_kf(t)\|_{L^q(0,T;L^p)}\Bigr\|_{\ell^r}.$$
\end{Definition}
Obviously,
$
\widetilde{L}^1_T(\dot{B}^s_{p,1})=L^1_T(\dot{B}^s_{p,1}).
$
In the sequel, we will constantly use the Bony's
decomposition from \cite{Bony} that
\beq\label{Bonydecom}
uv=T_uv+T_vu+R(u,v), \eeq with
$$T_uv=\sum_{j\in\Z}S_{j-1}u\Delta_jv, \quad R(u,v)=\sum_{j\in\Z}\Delta_ju \widetilde{\Delta}_{j}v,
\quad \widetilde{\Delta}_{j}v=\sum_{|j'-j|\le1}\Delta_{j'}v.$$

Let us conclude this section by collecting some useful lemmas.

\begin{Lemma}\cite{Che-book}\label{Lem:Bernstein}
Let $1\le p\le q\le+\infty$. Assume that $f\in L^p(\R^N)$, then for
any $\gamma\in(\N\cup\{0\})^N$, there exist constants $C_1$, $C_2$
independent of $f$, $j$ such that \beno &&{\rm supp}\hat f\subseteq
\{|\xi|\le A_02^{j}\}\Rightarrow \|\partial^\gamma f\|_q\le
C_12^{j{|\gamma|}+j N(\frac{1}{p}-\frac{1}{q})}\|f\|_{p},
\\
&&{\rm supp}\hat f\subseteq \{A_12^{j}\le|\xi|\le
A_22^{j}\}\Rightarrow \|f\|_{p}\le
C_22^{-j|\gamma|}\sup_{|\beta|=|\gamma|}\|\partial^\beta f\|_p.
\eeno
\end{Lemma}

\begin{Lemma}\cite{Dan-cpde}\label{Lem:poinc}
Let $1<p<\infty$, and $a\ge \bar{a}>0$ be a bounded continuous function.
Assume that $u\in L^p(\R^N)$ and $\textrm{supp}\,\hat u\subset \{\xi:R_1\le |\xi|\le R_2\}$.
Then there exists a constant $c$ depending only on $N$ and $R_2/R_1$ such that
\beno
c\bar{a}R_1^2\f{(p-1)} {p^2}\int_{\R^N}|u|^p\intd x\le -\int_{\R^N}{\rm div}(a\na u)|u|^{p-2}u \intd x.
\eeno
\end{Lemma}

\begin{Lemma}\label{Lem:binesti-1} Let $s>0$,
and $1\le p\le \infty$.
Assume that  $f, g\in \dot B^{s_1}_{p,1}\cap L^\infty$. Then there holds
\beno \|fg\|_{\dot B^{s}_{p,1}} \le C(\|f\|_{\dot
B^{s}_{p,1}}\|g\|_{{L}^{\infty}}+\|f\|_{{L}^{\infty}}\|g\|_{\dot B^{s}_{p,1}}).
\eeno
\end{Lemma}

\begin{Lemma}\label{Lem:binesti} Let $s_1, s_2\le \frac{N}{p},\, s_1+s_2>N\max (0,\frac2p-1)$,
and $1\le p,q,q_1,q_2\le \infty$ with $\f 1{q_1}+\f1{q_2}=\f1q$.
Assume that  $f\in \widetilde{L}^{q_1}_T(\dot B^{s_1}_{p,1})$ and
$g\in \widetilde{L}^{q_2}_T(\dot B^{s_2}_{p,1})$. Then there holds
\beno \|fg\|_{\widetilde{L}^{q}_T(\dot B^{s_1+s_2-
\frac{N}{p}}_{p,1})} \le C\|f\|_{\widetilde{L}^{q_1}_T(\dot
B^{s_1}_{p,1})}\|g\|_{\widetilde{L}^{q_2}_T(\dot B^{s_2}_{p,1})}.
\eeno
\end{Lemma}

\begin{Lemma}\label{Lem:binesti-end} Let $s_1\le \frac{N}{p}, s_2< \frac{N}{p},\, s_1+s_2\ge N\max(0,\frac 2p-1)$,
and $1\le p,q,q_1,q_2\le \infty$ with $\f 1{q_1}+\f1{q_2}=\f1q$.
Assume that  $f\in \widetilde{L}^{q_1}_T(\dot B^{s_1}_{p,1})$ and
$g\in \widetilde{L}^{q_2}_T(\dot B^{s_2}_{p,\infty})$. Then there
holds \beno \|fg\|_{\widetilde{L}^{q}_T(\dot B^{s_1+s_2-
\frac{N}{p}}_{p,\infty})} \le C\|f\|_{\widetilde{L}^{q_1}_T(\dot
B^{s_1}_{p,1})}\|g\|_{\widetilde{L}^{q_2}_T(\dot
B^{s_2}_{p,\infty})}. \eeno
\end{Lemma}

\begin{Lemma}\label{Lem:commu} Let $s\in (- N\min\big(\frac{1}{p},\frac{1}{p'}\big), \frac{N}{p}+1]$ and $1\le p,q,q_1,q_2\le \infty$ with $\f 1{q_1}+\f1{q_2}=\f1q$. Assume that $f\in \widetilde{L}^{q_1}_T(\dot B^{ \frac{N}{p}+1}_{p,1})$
and $g\in \widetilde{L}^{q_2}_T(\dot B^{s}_{p,1})$. Then  there
holds \beno \sum_j2^{j(s-1)}\|{\rm div}[\Delta_j,f]\na g\|_{L^q_T(L^p)}\le
C\|f\|_{\widetilde{L}^{q_1}_T(\dot B^{
\frac{N}{p}+1}_{p,1})}\|g\|_{\widetilde{L}^{q_2}_T(\dot
B^{s}_{p,1})}. \eeno
\end{Lemma}

\begin{Lemma}\label{Lem:nonesti}
Let $s>0$ and $1\le p,q\le \infty$. Assume that $F\in W^{[s]+3,\infty}_{loc}(\R)$ with  $F(0)=0$. Then
for any $f\in L^\infty_T(L^\infty)\cap \widetilde{L}^q_T(\dot B^s_{p,1})$, we have
\beno
\|F(f)\|_{\widetilde{L}^q_T(\dot B^s_{p,1})}\le C\bigl(1+\|f\|_{L^\infty_T(L^\infty)}\bigr)^{[s]+2}\|f\|_{\widetilde{L}^q_T(\dot B^s_{p,1})}.
\eeno
\end{Lemma}

Lemma \ref{Lem:binesti}-Lemma \ref{Lem:nonesti} can be easily proved
by using Bony's decomposition and Lemma \ref{Lem:Bernstein}, see
also \cite{Che-JDM, Dan-cpde} or Section 3 for similar results.

\begin{Remark}
Lemma \ref{Lem:binesti}-Lemma \ref{Lem:nonesti} still remain true
for the usual homogenous Besov spaces. For example, the estimate in
Lemma \ref{Lem:binesti} becomes \beno \|fg\|_{\dot B^{s_1+s_2-
\frac{N}{p}}_{p,1}} \le C\|f\|_{\dot B^{s_1}_{p,1}}\|g\|_{\dot
B^{s_2}_{p,1}}, \eeno with $p,s_1,s_2$ satisfying the conditions as
in Lemma \ref{Lem:binesti}.
\end{Remark}

\setcounter{equation}{0}
\section{Nonlinear estimates in the weighted Besov spaces}

Let us firstly introduce the weight function. Let $\{e_k(t)\}_{k\in\Z}$ be a sequence defined in $[0,+\infty)$
satisfying the following conditions:
\ben\label{assu:e_k}
e_k(t)\in [0,1], \quad e_k(t)\le e_{k'}(t)\quad \textrm{if}\quad k\le k'\quad \textrm{and} \quad
e_k(t)\sim e_{k'}(t)\quad \textrm{if}\quad k\sim k'.
\een
Then the weight function $\{\om_k(t)\}_{k\in\Z}$ is defined by
$$
\om_k(t)=\sum_{\ell\ge k}2^{k-\ell}e_{\ell}(t),\quad k\in \Z.
$$
It is easy to verify that for any $k\in \Z$,
\begin{equation} \label{weightprop}
\begin{split}
&\om_k(t)\le 2,\quad e_k(t)\le \om_k(t), \\
&\om_k(t)\le 2^{k-k'}\om_{k'}(t)\quad \textrm{if}\,\, k\ge k',\quad
\om_k(t)\le 3\om_{k'}(t)\quad \textrm{if}\,\, k\le k',
\\& \om_k(t)\sim \om_{k'}(t)\quad \textrm{if}\,\, k\sim k'.
\end{split}
\end{equation}

\begin{Definition}\label{Def:weightspacefuction}
Let $s\in\R$, $1\le p,r\le+\infty$, $0<T<+\infty$. The weighted
Besov space $\dot B^s_{p,r}(\om)$ is  defined by
$$\dot B^s_{p,r}(\om)=\{f\in {\cZ'}(
\R^{N}):\,\|f\|_{\dot
B^s_{p,r}(\om)}<+\infty\},$$ where
\beno
\|f\|_{\dot B^s_{p,r}(\om)}\eqdefa\bigl\|2^{ks}\om_k(T)
\|\Delta_kf\|_{p}\bigr\|_{\ell^r}. \eeno
\end{Definition}

\begin{Definition}\label{Def:weighttimespacefuction}
Let $s\in\R$, $1\le p,q\le+\infty$, $0<T<+\infty$. The weighted
function space $\widetilde{L}^q_T(\dot B^s_{p,1}(\om))$ is  defined
by
$$\widetilde{L}^q_T(\dot B^s_{p,1}(\om))=\{f\in L^q_T(\dot B^s_{p,1}(\om)):\,\|f\|_{\widetilde{L}^q_T(\dot
B^s_{p,1}(\om))}<+\infty\},$$ where \beno
\|f\|_{\widetilde{L}^q_T(\dot B^s_{p,1}(\om))}\eqdefa\sum_{k\in\Z}2^{ks}\om_k(T)
\bigg(\int_0^T\|\Delta_kf(t)\|^q_{p}dt\bigg)
^{\frac{1}{q}}. \eeno
\end{Definition}

\begin{Remark}\label{rem:smallnorm}
If $e_k(t)$ is continuous on $[0,+\infty)$ and $e_k(0)=0$ for
$k\in\Z$, $f\in \widetilde{L}^\infty_T(\dot B^{s}_{p,1})$, then for
any $\varepsilon>0$, there exists a $\widetilde{T}\in (0,T]$ such
that \beno \|f\|_{\widetilde{L}^\infty_{\widetilde{T}}(\dot
B^s_{p,1}(\om))}\le \varepsilon. \eeno Indeed, due to $f\in
\widetilde{L}^\infty_T(\dot B^{s}_{p,1})$ and $\om_k(T)\le 2$, there
exists $N_1\in\N$ such that \beno \sum_{|k|\ge
N_1+1}2^{ks}\om_k(T)\|\Delta_k f\|_{L^\infty_T(L^p)}\le
\varepsilon/3,\quad\; \sum_{|k|\le N_1}2^{ks}\sum_{\ell\ge
k+N_1+1}2^{k-\ell}e_{\ell}(T)\|\Delta_k f\|_{L^\infty_T(L^p)}\le
\varepsilon/3. \eeno Thus, we have \beno
\|f\|_{\widetilde{L}^\infty_{\widetilde{T}}(\dot
B^s_{p,1}(\om))}&\le& 2\varepsilon/3
+\sum_{|k|\le N_1}2^{ks}\sum_{k\le\ell\le k+N_1}2^{k-\ell}e_{\ell}(\widetilde{T})\|\Delta_k f\|_{L^\infty_{\widetilde{T}}(L^p)}\\
&\le& 2\varepsilon/3
+2e_{2N_1}(\widetilde{T})\sum_{|k|\le N_1}2^{ks}\|\Delta_k f\|_{L^\infty_{\widetilde{T}}(L^p)}\\
&\le& \varepsilon, \eeno if $\widetilde{T}\in (0,T]$ is chosen such
that \beno 2e_{2N_1}(\widetilde{T})\sum_{|k|\le N_1}2^{ks}\|\Delta_k
f\|_{L^\infty_{\widetilde{T}}(L^p)}\le \varepsilon/3. \eeno
\end{Remark}

Next, we present some estimates in the weighted Besov spaces.

\begin{Lemma}\label{Lem:bonyweighbesov} Let $1\le p\le \infty$.
Assume that  $f\in \dot B^{s_1}_{p,1}(\omega), g\in \dot
B^{s_2}_{p,1}$. Then there hold

\no (a)\, if $s_2\le \frac{N}{p}$, we have \beno \|T_gf\|_{\dot
B^{s_1+s_2- \frac{N}{p}}_{p,1}(\omega)} \le C\|f\|_{\dot
B^{s_1}_{p,1}(\omega)}\|g\|_{\dot B^{s_2}_{p,1}}; \eeno

\no (b)\, if $s_1\le \frac{N}{p}-1$, we have \beno \|T_fg\|_{\dot
B^{s_1+s_2- \frac{N}{p}}_{p,1}(\omega)} \le C\|f\|_{\dot
B^{s_1}_{p,1}(\omega)}\|g\|_{\dot B^{s_2}_{p,1}}; \eeno

\no (c)\, if $s_1+s_2>N\max(0,\frac 2p-1)$, we have \beno
\|R(f,g)\|_{\dot B^{s_1+s_2- \frac{N}{p}}_{p,1}(\omega)} \le
C\|f\|_{\dot B^{s_1}_{p,1}(\omega)}\|g\|_{\dot B^{s_2}_{p,1}}. \eeno

\end{Lemma}

\no{\bf Proof.}\, Due to (\ref{orth}), we have
\beno
\Delta_j(T_gf)=\sum_{|j'-j|\le 4}\Delta_j(S_{j'-1}g\Delta_{j'}f),
\eeno
then we get by Lemma \ref{Lem:Bernstein} and (\ref{weightprop}) that
\beno
\|T_gf\|_{\dot B^{s_1+s_2- \frac{N}{p}}_{p,1}(\omega)}&=&\sum_{j}2^{j(s_1+s_2- \frac{N}{p})}\om_j(T)\|\Delta_j(T_gf)\|_{p}\\
&\le& C\sum_{j}2^{j(s_1+s_2- \frac{N}{p})}\om_{j}(T)\|S_{j-1}g\|_{\infty}\|\Delta_{j}f\|_{p}\\
&\le& C\|f\|_{\dot B^{s_1}_{p,1}(\omega)}\|g\|_{\dot B^{s_2}_{p,1}},
\eeno where we used in the last inequality \beno
\|S_{j-1}g\|_{\infty} \le C\sum_{\ell\le j-2}2^{\ell
 \frac{N}{p}}\|\Delta_{\ell}g\|_{p}\le C2^{j(-s_2+ \frac{N}{p})}\|g\|_{\dot
B^{s_2}_{p,1}}. \eeno This proves (a). We next prove (b). Similarly,
we have \beno
\|T_fg\|_{\dot B^{s_1+s_2- \frac{N}{p}}_{p,1}(\omega)}&=&\sum_{j}2^{j(s_1+s_2- \frac{N}{p})}\om_j(T)\|\Delta_j(T_fg)\|_{p}\\
&\le& C\sum_{j}2^{j(s_1+s_2-
\frac{N}{p})}\om_{j}(T)\|S_{j-1}f\|_{\infty}\|\Delta_{j}g\|_{p},
\eeno and by Lemma \ref{Lem:Bernstein} and (\ref{weightprop}), we
have \beno
\om_{j}(T)\|S_{j-1}f\|_{\infty}&\le& C2^{j}\sum_{\ell\le j-2}2^{\ell ( \frac{N}{p}-1)}\om_{\ell}(T)\|\Delta_{\ell}f\|_{p}\\
&\le& C2^{j( \frac{N}{p}-s_1)}\|f\|_{\dot B^{s_1}_{p,1}(\omega)},
\eeno which lead to (b). Now we prove (c). Notice that \beno
\Delta_j(R(f,g))=\sum_{j'\ge
j-3}\Delta_j(\Delta_{j'}f\widetilde{\Delta}_{j'}g), \eeno then we get by Lemma
\ref{Lem:Bernstein} that if $p\ge2$\beno
&&\|R(f,g)\|_{\dot B^{s_1+s_2- \frac{N}{p}}_{p,1}(\omega)}\\
&&\le C\sum_{j}\sum_{j'\ge j-3}2^{j(s_1+s_2)}\om_j(T)\|\Delta_{j'}f\|_{p}\|\widetilde{\Delta}_{j'}g\|_{p}\\
&&\le C\sum_{j}\sum_{j'\ge j-3}\sum_{\ell\ge j}2^{j-\ell}e_\ell(T)2^{j(s_1+s_2)}\|\Delta_{j'}f\|_{p}\|\widetilde{\Delta}_{j'}g\|_{p}\\
&&= \sum_{j}\sum_{j'\ge j-3}\sum_{\ell\ge j,j'}\square+\sum_{j}\sum_{j'\ge j-3}\sum_{j'\ge \ell\ge j}\square\\
&&\triangleq I+II. \eeno

For $II$, using the fact that $e_\ell(T)\le e_{j'}(T)\le
\om_{j'}(T)$ if $\ell\le j'$, we get \beno
II&\le& C\sum_{j}\sum_{j'\ge j-3}\om_{j'}(T)2^{j(s_1+s_2)}\|\Delta_{j'}f\|_{p}\|\widetilde{\Delta}_{j'}g\|_{p}\\
&\le& C\sum_{j}\sum_{j'\ge j-3}\om_{j'}(T)2^{(j-j')(s_1+s_2)}2^{j's_1}\|\Delta_{j'}f\|_{p}\|g\|_{\dot B^{s_2}_{p,1}}\\
&\le& C\|f\|_{\dot B^{s_1}_{p,1}(\omega)}\|g\|_{\dot B^{s_2}_{p,1}},
\eeno and for $I$, using the fact that \beno \sum_{\ell\ge
j,j'}2^{j-\ell}e_\ell(T)\le 2^{j-j'}\sum_{\ell\ge
j'}2^{j'-\ell}e_\ell(T)=2^{j-j'}w_{j'}(T), \eeno we obtain \beno
I&\le& C\sum_{j}\sum_{j'\ge j-3}\om_{j'}(T)2^{j(s_1+s_2)}2^{j-j'}\|\Delta_{j'}f\|_{p}\|\widetilde{\Delta}_{j'}g\|_{p}\\
&\le& C\sum_{j}\sum_{j'\ge j-3}\om_{j'}(T)2^{(j-j')(s_1+s_2+1)}2^{j's_1}\|\Delta_{j'}f\|_{p}\|g\|_{\dot B^{s_2}_{p,1}}\\
&\le& C\|f\|_{\dot B^{s_1}_{p,1}(\omega)}\|g\|_{\dot B^{s_2}_{p,1}}.
\eeno  If $p<2$, we get by Lemma \ref{Lem:Bernstein} that \beno
&&\|R(f,g)\|_{\dot B^{s_1+s_2- \frac{N}{p}}_{p,1}(\omega)}\\
&&\le C\sum_{j}\sum_{j'\ge j-3}2^{j(s_1+s_2-N(\frac2p-1))}\om_j(T)\|\Delta_{j'}f\|_{p}\|\widetilde{\Delta}_{j'}g\|_{p'}\\
&&\le C\sum_{j}\sum_{j'\ge j-3}\sum_{\ell\ge
j}2^{j-\ell}e_\ell(T)2^{j(s_1+s_2-N(\frac2p-1))}\|\Delta_{j'}f\|_{p}
\|\widetilde{\Delta}_{j'}g\|_{p}2^{N(\frac2p-1)j'}. \eeno Then
treating it as in the case of $p\ge2$, we  obtain the same
inequality  for $s_1+s_2> N(\frac2p-1)$. This proves (c). \ef

We have a similar result in the weighted Besov spaces with the time.

\begin{Lemma}\label{Lem:bonyweighttimespace} Let $1\le p,q,q_1,q_2\le \infty$ with $\f 1{q_1}+\f1{q_2}=\f1q$.
Assume that  $f\in \widetilde{L}^{q_1}_T(\dot
B^{s_1}_{p,1}(\omega)), g\in \widetilde{L}^{q_2}_T(\dot
B^{s_2}_{p,1})$. Then there hold

\no (a)\, if $s_2\le \frac{N}{p}$, we have \beno
\|T_gf\|_{\widetilde{L}^{q}_T(\dot B^{s_1+s_2-
\frac{N}{p}}_{p,1}(\omega))} \le C\|f\|_{\widetilde{L}^{q_1}_T(\dot
B^{s_1}_{p,1}(\omega))}\|g\|_{\widetilde{L}^{q_2}_T(\dot
B^{s_2}_{p,1})}; \eeno

\no (b)\, if $s_1\le \frac{N}{p}-1$, we have \beno
\|T_fg\|_{\widetilde{L}^{q}_T(\dot B^{s_1+s_2-
\frac{N}{p}}_{p,1}(\omega))} \le C\|f\|_{\widetilde{L}^{q_1}_T(\dot
B^{s_1}_{p,1}(\omega))}\|g\|_{\widetilde{L}^{q_2}_T(\dot
B^{s_2}_{p,1})}; \eeno

\no (c)\, if $s_1+s_2>N\max(0,\frac 2p-1)$, we have \beno
\|R(f,g)\|_{\widetilde{L}^{q}_T(\dot B^{s_1+s_2-
\frac{N}{p}}_{p,1}(\omega))} \le C\|f\|_{\widetilde{L}^{q_1}_T(\dot
B^{s_1}_{p,1}(\omega))}\|g\|_{\widetilde{L}^{q_2}_T(\dot
B^{s_2}_{p,1})}. \eeno

\end{Lemma}

The following proposition is a direct consequence of Lemma \ref{Lem:bonyweighttimespace}.

\begin{Proposition}\label{Prop:binweight} Let $s_1\le \frac{N}{p}-1, s_2\le \frac{N}{p}, s_1+s_2>N\max(0,\frac 2p-1)$, and $1\le p,q,q_1,q_2\le \infty$
with $\f 1{q_1}+\f1{q_2}=\f1q$. Assume that  $f\in
\widetilde{L}^{q_1}_T(\dot B^{s_1}_{p,1}(\omega))$ and $g\in
\widetilde{L}^{q_2}_T(\dot B^{s_2}_{p,1})$. Then there holds \beno
\|fg\|_{\widetilde{L}^{q}_T(\dot B^{s_1+s_2-
\frac{N}{p}}_{p,1}(\omega))} \le C\|f\|_{\widetilde{L}^{q_1}_T(\dot
B^{s_1}_{p,1}(\omega))}\|g\|_{\widetilde{L}^{q_2}_T(\dot
B^{s_2}_{p,1})}. \eeno
\end{Proposition}

From the proof of Lemma \ref{Lem:bonyweighbesov}, we can also obtain

\begin{Proposition}\label{Prop:binweight-end} Let $s_1\le \frac{N}{p}-1, s_2< \frac{N}{p}, s_1+s_2\ge N\max(0,\frac 2p-1)$, and $1\le p,q,q_1,q_2\le \infty$
with $\f 1{q_1}+\f1{q_2}=\f1q$. Assume that  $f\in
\widetilde{L}^{q_1}_T(\dot B^{s_1}_{p,1}(\omega))$ and $g\in
\widetilde{L}^{q_2}_T(\dot B^{s_2}_{p,\infty})$. Then there holds
\beno \|fg\|_{\widetilde{L}^{q}_T(\dot B^{s_1+s_2-
\frac{N}{p}}_{p,\infty}(\om))} \le
C\|f\|_{\widetilde{L}^{q_1}_T(\dot
B^{s_1}_{p,1}(\omega))}\|g\|_{\widetilde{L}^{q_2}_T(\dot
B^{s_2}_{p,\infty})}. \eeno
\end{Proposition}

\begin{Proposition}\label{Prop:nonweight}
Let $s>0$ and $1\le p,q\le \infty$. Assume that $F\in W^{[s]+3,\infty}_{loc}(\R)$ with  $F(0)=0$. Then
for any $f\in L^\infty_T(L^\infty)\cap \widetilde{L}^q_T(\dot B^s_{p,1}(\om))$, we have
\beno
\|F(f)\|_{\widetilde{L}^q_T(\dot B^s_{p,1}(\om))}\le C(1+\|f\|_{L^\infty_T(L^\infty)})^{[s]+2}\|f\|_{\widetilde{L}^q_T(\dot B^s_{p,1}(\om))}.\label{5.12}
\eeno
\end{Proposition}

\no{\bf Proof.}\, We decompose $F(f)$ as
\begin{align}
F(f)=\sum_{j'\in\Z}F(S_{j'+1}f)-F(S_{j'}f)&=\sum_{j'\in\Z}\Delta_{j'}f\int_0^1F'(S_{j'}f+\tau\Delta_{j'}f)d\tau\nonumber\\
&\triangleq\sum_{j'\in\Z}\Delta_{j'}f\, m_{j'}(f),\nonumber
\end{align}
where $m_{j'}(f)=\int_0^1F'(S_{j'}f+\tau\Delta_{j'}f)d\tau$. Furthermore, we write
\begin{align}
\Delta_jF(f)=\sum_{j'<j}\Delta_j(\Delta_{j'}f\,m_{j'}(f))+\sum_{j'\ge j}\Delta_j(\Delta_{j'}f\, m_{j'}(f))\triangleq I+II.\nonumber
\end{align}
By Lemma \ref{Lem:Bernstein}, we have
\begin{align}\label{eq:I-localestimate}
\|I\|_{L^q_T(L^p)}&\le\sum_{j'<j}\|\Delta_j(\Delta_{j'}f\,m_{j'}(f))\|_{L^q_T(L^p)}\nonumber\\
&\le \sum_{j'<j}2^{-j|\al|}\sup_{|\gamma|=|\al|}\|D^\gamma\Delta_j(\Delta_{j'}f\,m_{j'}(f))\|_{L^q_T(L^p)},
\end{align}
with $\al$ to be determined later. Notice that for $|\gamma|\ge0$, we have
$$\|D^\gamma m_{j'}(f)\|_\infty\le C2^{j'|\gamma|}(1+\|f\|_\infty)^{|\gamma|}\|F'\|_{W^{|\gamma|,\infty}},$$
from which and (\ref{eq:I-localestimate}), it follows that
\begin{align}
2^{js}\|I\|_{L^q_T(L^p)}\le C2^{j(s-|\al|)}\sum_{j'<j}2^{j'|\al|}
\|\Delta_{j'}f\|_{L^q_T(L^p)}(1+\|f\|_{L^\infty_T(L^\infty)})^{|\al|}\|F'\|_{W^{|\al|,\infty}},\nonumber
\end{align}
thus, if we take $|\al|=[s]+2$, we get by (\ref{weightprop}) that
\begin{align}\label{eq:I-weightedestimate}
&\sum_{j}\omega_j(T)2^{js}\|I\|_{L^q_T(L^p)}\nonumber\\
&\le C\sum_{j'}2^{j's}\omega_{j'}(T)\|\Delta_{j'}f\|_{L^q_T(L^p)}
\sum_{j>j'}2^{(j-j')(s-|\al|+1)}(1+\|f\|_{L^\infty_T(L^\infty)})^{|\al|}\|F'\|_{W^{|\al|,\infty}}
\nonumber\\&\le C(1+\|f\|_{L^\infty_T(L^\infty)})^{[s]+2}\|F'\|_{W^{[s]+2,\infty}}\|f\|_{\widetilde{L}^q_T(\dot B^s_{p,1}(\om))}.
\end{align}

Now, let us turn to the estimate of $II$. We get by Lemma \ref{Lem:Bernstein} that
\begin{align}
\|II\|_{L^q_T(L^p)}&\le C\sum_{j'\ge j}\|\Delta_{j'}f\|_{L^q_T(L^p)}.\nonumber
\end{align}
Then we write
\begin{align}
\sum_{j}\omega_j(T)2^{js}\|II\|_{L^q_T(L^p)}
&\le C\sum_{j}2^{js}\sum_{j'\ge j}\|\Delta_{j'}f\|_{L^q_T(L^p)}
\sum_{j'\ge \ell \ge j}2^{j-\ell}e_{\ell}(T)\nonumber\\
&\quad+C\sum_{j}2^{js}\sum_{j'\ge j}\|\Delta_{j'}f\|_{L^q_T(L^p)}
\sum_{\ell \ge j, j'}2^{j-\ell}e_{\ell}(T),\nonumber
\end{align}
from which and a similar argument of (c) in Lemma \ref{Lem:bonyweighbesov}, we infer that
\beno
\sum_{j}\omega_j(T)2^{js}\|II\|_{L^q_T(L^p)}\le C\|f\|_{\widetilde{L}^q_T(\dot B^s_{p,1}(\om))},
\eeno
from which and (\ref{eq:I-weightedestimate}), we conclude the proof of Proposition \ref{Prop:nonweight}.\ef

\setcounter{equation}{0}
\section{Estimates of the linear transport equation}

In this section, we study  the linear transport equation
\begin{align}\label{eq:linertrans}
\bigg\{\begin{aligned}
&\partial_t f+v\cdot \nabla f =g,\\
&f(0,x)=f_0.
\end{aligned}
\bigg.\end{align}

\begin{Proposition}\cite{Dan-NDEA}\label{Prop:transport}
Let $s\in (-N\min(\frac1p,\frac1{p'}), 1+\frac Np)$, $1\le
p,r\le+\infty$, and $s=1+\frac Np$, if $r=1$. Let $v$ be a vector field such that $\nabla v\in
L^1_T(\dot{B}^{\frac{N}{p}}_{p,r}\cap L^\infty)$. Assume that
$f_0\in \dot{B}^{s}_{p,r},$ $g\in L^1_T(\dot{B}^{s}_{p,r})$ and $f$
is the solution of (\ref{eq:linertrans}). Then there holds for
$t\in[0,T]$, \beno
\|f\|_{\widetilde{L}^\infty_t(\dot{B}^{s}_{p,r})}\le e^{CV(t)}\bigg(
\|f_0\|_{\dot{B}^{s}_{p,r}}+\int_0^t
e^{-CV(\tau)}\|g(\tau)\|_{\dot{B}^{s}_{p,r}}d\tau\bigg), \eeno where
$V(t)\eqdefa\int_0^t\|\nabla
v(\tau)\|_{\dot{B}^{\frac{N}{p}}_{p,r}\cap L^\infty}d\tau.$ If
$r<+\infty$, then $f$ belongs to $C([0,T]; \dot{B}^s_{p,r})$.
\end{Proposition}

\begin{Proposition}\label{Prop:transportweight}
Let $p\in [1,+\infty]$ and $s\in (-N\min(\frac1p,\frac1{p'}),\frac
Np]$. Let $v$ be a vector field such that $\nabla v\in
L^1_T(\dot{B}^{ \frac{N}{p}}_{p,1})$. Assume that $f_0\in
\dot{B}^{s}_{p,1},$ $g\in L^1_T(\dot{B}^{s}_{p,1})$ and $f$ is the
solution of (\ref{eq:linertrans}). Then there holds for $t\in
[0,T]$, \beno \|f\|_{\widetilde{L}^\infty_t(\dot
B^s_{p,1}(\omega))}\le e^{CV(t)}\bigg(
\|f_0\|_{\dot{B}^{s}_{p,1}(\om)}+\int_0^t
e^{-CV(\tau)}\|g(\tau)\|_{\dot{B}^{s}_{p,1}(\om)}d\tau\bigg), \eeno
where $V(t)\eqdefa\int_0^t\|\nabla
v(\tau)\|_{\dot{B}^{\frac{N}{p}}_{p,1}}d\tau.$
\end{Proposition}

\no{\bf Proof.}\, Applying the operator $\Delta_j$ to the transport equation, we obtain
\ben\label{eq:loctrans}
\p_t\Delta_jf+v\cdot\na \Delta_jf=\Delta_jg+[v,\Delta_j]\cdot \na f.
\een
Assume that $p<+\infty$. Multiplying both sides of (\ref{eq:loctrans}) by $|\Delta_jf|^{p-2}\Delta_jf$,
we get by integrating by parts over $\R^N$ for the resulting equation that
\beno
\f1p\f{d}{dt}\|\Delta_jf\|_p^p-\f1p\int_{\R^N}|\Delta_jf|^{p}\dv vdx\le \bigl(\|\Delta_jg\|_p+\|[v,\Delta_j]\cdot \na f\|_p\bigr)\|\Delta_jf\|_p^{p-1},
\eeno
then we have
\beno
\|\Delta_jf(t)\|_p\le \|\Delta_jf_0\|_p+\int_0^t\bigl(\|\Delta_jg\|_p+\|[v,\Delta_j]\cdot \na f\|_p+\f1p\|\dv v\|_\infty\|\Delta_jf\|_p\bigr)d\tau,
\eeno
from which, it follows that
\beno
\|f\|_{\widetilde{L}^\infty_t(\dot B^s_{p,1}(\omega))}&\le&
\|f_0\|_{\dot B^s_{p,1}(\om)}+C\int_0^t\|\dv v(\tau)\|_{\dot B^{ \frac{N}{p}}_{p,1}}\|f(\tau)\|_{\widetilde{L}^\infty_\tau(\dot B^s_{p,1}(\omega))}d\tau\nonumber\\
&&+\int_0^t\|g(\tau)\|_{\dot
B^s_{p,1}(\omega)}d\tau+\int_0^t\sum_{j}\om_j(T)2^{js}\|[v,\Delta_j]\cdot
\na f(\tau)\|_pd\tau, \eeno from which and Lemma
\ref{Lem:commweight}, we infer that \beno
\|f\|_{\widetilde{L}^\infty_t(\dot B^s_{p,1}(\omega))}\le
\|f_0\|_{\dot B^s_{p,1}(\om)}+C\int_0^t\|v(\tau)\|_{\dot B^{
\frac{N}{p}+1}_{p,1}}\|f(\tau)\|_{\widetilde{L}^\infty_\tau(\dot
B^s_{p,1}(\omega))}d\tau +\int_0^t\|g(\tau)\|_{\dot
B^s_{p,1}(\omega)}d\tau. \eeno Then Gronwall's lemma applied implies
the desired inequality.\ef

\begin{Lemma}\label{Lem:commweight} Let $p\in [1,\infty], s\in (-N\min(\frac1p,\frac1{p'}), \frac Np]$.
Assume that $v\in \dot B^{ \frac{N}{p}+1}_{p,1}$ and $f\in \dot
B^s_{p,1}(\om)$. Then there holds \beno
\sum_{j}\om_j(T)2^{js}\|[v,\Delta_j]\cdot \na f\|_p\le C\|v\|_{\dot
B^{ \frac{N}{p}+1}_{p,1}}\|f\|_{\dot B^s_{p,1}(\om)}. \eeno
\end{Lemma}

\no{\bf Proof.}\, Using the Bony's decomposition, we write
\beno
[v,\Delta_j]\cdot \na f&=&[T_{v^k},\Delta_j]\p_kf+T_{\p_k\Delta_jf}v^k+R(v^k,\p_k\Delta_jf)\\
&&-\Delta_j(T_{\p_kf}v^k)-\Delta_jR(v^k,\p_kf). \eeno Using Lemma
\ref{Lem:bonyweighbesov} with $s_1=s-1$ and $s_2= \frac{N}{p}+1$, we
get \beno
&&\sum_{j}\om_j(T)2^{js}\|\Delta_j(T_{\p_kf}v^k)\|_p\le C\|v\|_{\dot B^{ \frac{N}{p}+1}_{p,1}}\|f\|_{\dot B^s_{p,1}(\om)},\\
&&\sum_{j}\om_j(T)2^{js}\|\Delta_jR(v^k,\p_kf)\|_p\le C\|v\|_{\dot
B^{ \frac{N}{p}+1}_{p,1}}\|f\|_{\dot B^s_{p,1}(\om)}. \eeno Notice
that \beno T_{\p_k\Delta_jf}'v^k\triangleq
T_{\p_k\Delta_jf}v^k+R(v^k,\p_k\Delta_jf)=\sum_{j'\ge
j-2}S_{j'+2}\Delta_j\p_kf\Delta_{j'}v^k, \eeno then we get by Lemma
\ref{Lem:Bernstein} that \beno
\sum_{j}\om_j(T)2^{js}\|T_{\p_k\Delta_jf}'v^k\|_p&\le& C\sum_{j}\om_j(T)2^{js}\|\Delta_j\na f\|_\infty\sum_{j'\ge j-2}\|\Delta_{j'}v^k\|_p\\
&\le& C\sum_{j}\om_j(T)2^{j(s+1+ \frac{N}{p})}\|\Delta_jf\|_p\sum_{j'\ge j}\|\Delta_{j'}v^k\|_p\\
&\le& C\|v\|_{\dot B^{ \frac{N}{p}+1}_{p,1}}\|f\|_{\dot
B^s_{p,1}(\om)}. \eeno Now, we turn to estimate $[T_{v^k},
\Delta_j]\p_k f$.  Set $h(x)=(\cF^{-1}\varphi)(x)$, we get by using
Taylor's formula that \beno
[T_{v^k}, \Delta_j]\pa_k f&=&\sum_{|j'-j|\le4}[S_{j'-1}v^k, \Delta_j]\pa_k \Delta_{j'}f\\
&=&\sum_{|j'-j|\le4}2^{Nj}\int_{\R^N}h(2^j(x-y))(S_{j'-1}v^k(x)-S_{j'-1}v^k(y))\pa_k \Delta_{j'}f(y)dy\\
&=&\sum_{|j'-j|\le4}2^{(N+1)j}\int_{\R^N}\int_0^1y\cdot\na
S_{j'-1}v^k(x-\tau
y)d\tau\pa_kh(2^jy)\Delta_{j'}f(x-y)dy\nonumber\\&&\qquad\quad+
2^{Nj}\int_{\R^N}h(2^j(x-y))\pa_kS_{j'-1}v^k(y)\Delta_{j'}f(y)dy ,
\eeno from which and the Minkowski inequality, we infer that \beno
\sum_j\om_j(T)2^{js}\|[T_{v^k}, \Delta_j]\pa_k f\|_p &\le&
C\sum_j\om_j(T)2^{js}\sum_{|j'-j|\le4}\|\na
S_{j'-1}v\|_\infty\|\Delta_{j'}f\|_p\\
 &\le& C\|v\|_{\dot B^{ \frac{N}{p}+1}_{p,1}}\|f\|_{\dot B^s_{p,1}(\om)}.
\eeno

Summing up all the above estimates, we conclude the proof of Lemma \ref{Lem:commweight}.\ef

\setcounter{equation}{0}
\section{Estimates of the linearized momentum equation}

In this section, we study the linearized momentum equation
\ben\label{eq:linearmomen} \left\{
\begin{array}{ll}
\p_tu-\dv(\overline{\mu}\na u)
-\na((\overline{\lambda}+\overline{\mu})\dv\,u)=G, \\
u|_{t=0}=u_0.
\end{array}
\right. \een
In what follows, we assume that the viscosity coefficients $\overline{\lambda}(\rho)$ and $\overline{\mu}(\rho)$ depend smoothly on the function $\rho$
and there exists a positive constant $c_1$ such that
\beno
\overline{\mu}\ge c_1,\quad \overline{\lambda}+2\overline{\mu}\ge c_1.
\eeno
Fix a positive constant $c$ to be chosen later. In this section, the weighted function $\om_k(t)$ is given by
$$
\om_k(t)=\sum_{\ell\ge k}2^{k-\ell}e_\ell(t),
$$
with $e_\ell(t)=(1-e^{-c2^{2\ell}t})^\f12$. It is easy to verify that the function $e_\ell(t)$ satisfies (\ref{assu:e_k}).

\begin{Proposition}\label{Prop:momenequ} Let $q\in [1,\infty]$.
Assume that $G\in L^1_T(\dot B^{s-1}_{p,1}), u_0\in \dot
B^{s-1}_{p,1}$, and $\rho-\underline{\rho}\in L^\infty_T(\dot B^{
\frac{N}{p}}_{p,1})$. Let $u$ be a solution of
(\ref{eq:linearmomen}). Then there hold

\no (a)\, If $p\in (1,N],$ $s\in (-N\min(\frac1p,\frac1{p'})+1,\frac Np]$, we
have \beno \|u\|_{\widetilde{L}^q_T(\dot B^{s-1+2/q}_{p,1})} \le
C\Bigl(\|u_0\|_{\dot
B^{s-1}_{p,1}}+\|G(\tau)\|_{\widetilde{L}^1_T(\dot B^{s-1}_{p,1})}
+A(T)\|\rho-\underline{\rho}\|_{\widetilde{L}^\infty_T(\dot B^{
\frac{N}{p}}_{p,1})}\|u\|_{\widetilde{L}^1_T(\dot
B^{s+1}_{p,1})}\Bigr); \eeno In addition, if
$\rho-\underline{\rho}\in L^\infty_T(\dot B^{
\frac{N}{p}+1}_{p,1})$, $p\in(1,\infty)$, $s\in (-N\min(\frac1p,\frac1{p'})+1,\frac Np+1]$, then \beno \|u\|_{\widetilde{L}^q_T(\dot
B^{s-1+2/q}_{p,1})} \le C\Bigl(\|u_0\|_{\dot
B^{s-1}_{p,1}}+\|G(\tau)\|_{\widetilde{L}^1_T(\dot B^{s-1}_{p,1})}
+A(T)\|\rho-\underline{\rho}\|_{\widetilde{L}^\infty_T(\dot B^{
\frac{N}{p}+1}_{p,1})}\|u\|_{\widetilde{L}^1_T(\dot
B^{s}_{p,1})}\Bigr); \eeno

\no (b)\,If  $p\in (1,N],$  $s\in (-N\min(\frac1p,\frac1{p'})+1,\frac Np]$, we have
\beno &&\|u\|_{\widetilde{L}^1_T(\dot
B^{s+1}_{p,1})}+\|u\|_{\widetilde{L}^2_T(\dot B^{s}_{p,1})}\\&& \le
C\Bigl(\|u_0\|_{\dot B^{s-1}_{p,1}(\om))}
+\|G(\tau)\|_{\widetilde{L}^1_T(\dot B^{s-1}_{p,1}(\om))}
+A(T)\|\rho-\underline{\rho}\|_{\widetilde{L}^\infty_T(\dot B^{
\frac{N}{p}}_{p,1}(\om))}\|u\|_{\widetilde{L}^1_T(\dot
B^{s+1}_{p,1})}\Bigr). \eeno Here
$A(T)\eqdefa\bigl(1+\|\rho\|_{L^\infty_T(L^\infty)}\bigr)^{[
\frac{N}{p}]+2}$.

\end{Proposition}

\no{\bf Proof.}\, Set $d=\dv u$ and $w=\textrm{curl} u$. From (\ref{eq:linearmomen}), we find that $(d,w)$ satisfies
\ben\label{eq:dv-curl} \left\{
\begin{array}{ll}
\p_td-\dv(\overline{\nu}\na d)=\dv G+F_1, \\
\p_tw-\dv(\overline{\mu}\na w)=\textrm{curl} G+F_2,\\
(d,w)|_{t=0}=(\dv u_0,\textrm{ curl} u_0)\triangleq (d_0,w_0),
\end{array}
\right. \een
where $\overline{\nu}=\overline{\lambda}+2\overline{\mu}$ and
\beno
&&F_1=\dv(\na \overline{\mu}\cdot\na u)+\dv(\na (\overline{\lambda}+\overline{\mu})d),\\
&&F_2^{i,j}=\dv(\p_j\overline{\mu}\na u^i-\p_i\overline{\mu}\na u^j),\quad i,j=1,\cdots,N.
\eeno

Applying the operator $\Delta_j$ to (\ref{eq:dv-curl}), we obtain
\beno \left\{
\begin{array}{ll}
\p_t\Delta_jd-\dv(\overline{\nu}\na \Delta_jd)=\dv \Delta_jG+\Delta_jF_1+\dv[\Delta_j,\overline{\nu}]\na d, \\
\p_t\Delta_jw-\dv(\overline{\mu}\na \Delta_jw)=\textrm{curl} \Delta_jG+\Delta_jF_2
+\dv[\Delta_j,\overline{\mu}]\na w.
\end{array}
\right. \eeno
Multiplying the first equation by $|\Delta_jd|^{p-2}\Delta_jd$, we get by integrating over $\R^N$ that
\beno
&&\f1p\f \intd{\intd t}\|\Delta_jd\|_p^p-\int_{\R^N}\dv(\overline{\nu}\na \Delta_jd)|\Delta_jd|^{p-2}\Delta_jd \intd x\\
&&\quad=\int_{\R^N}\bigl(\dv
\Delta_jG+\Delta_jF_1+\dv[\Delta_j,\overline{\nu}]\na
d\bigr)|\Delta_jd|^{p-2}\Delta_jd \intd x \eeno Lemma
\ref{Lem:poinc} ensures there exists a positive constant $c_p$
depending on $c_0, p,N$ such that \beno \f1p\f \intd{\intd
t}\|\Delta_jd\|_p^p+c_p2^{2j}\|\Delta_jd\|_p^p\le
\|\Delta_jd\|_p^{p-1} \bigl(\|\dv
\Delta_jG\|_p+\|\Delta_jF_1\|_p+\|\dv[\Delta_j,\overline{\nu}]\na
d\|_p\bigr). \eeno Thus, we have \beno \f \intd{\intd
t}\|\Delta_jd\|_p+c_p2^{2j}\|\Delta_jd\|_p\le\|\dv
\Delta_jG\|_p+\|\Delta_jF_1\|_p+\|\dv[\Delta_j,\overline{\nu}]\na
d\|_p, \eeno which implies that \beno \|\Delta_jd(t)\|_p\le
e^{-c_p2^{2j}t}\|\Delta_jd_0\|_p+\int_0^te^{-c_p2^{2j}(t-\tau)}\bigl(\|\dv
\Delta_jG\|_p+\|\Delta_jF_1\|_p+\|\dv[\Delta_j,\overline{\nu}]\na
d\|_p\bigr)\intd\tau. \eeno Similarly, we can obtain \beno
\|\Delta_jw(t)\|_p\le
e^{-c_p2^{2j}t}\|\Delta_jw_0\|_p+\int_0^te^{-c_p2^{2j}(t-\tau)}\bigl(\|\textrm{curl}\Delta_jG\|_p+\|\Delta_jF_2\|_p+\|\dv[\Delta_j,\overline{\mu}]\na
w\|_p\bigr)\intd\tau. \eeno From the above two inequalities, we
infer that for any $q\in [1,\infty]$ and $t\in [0,T]$,
\ben\label{eq:localcurl-dvest}
&&\|\Delta_jd(t)\|_{L^q_t(L^p)}+\|\Delta_jw(t)\|_{L^q_t(L^p)}\nonumber\\
&&\le C2^{-2j/q}c_j(T)^\f1q(\|\Delta_jd_0\|_p+\|\Delta_jw_0\|_p)\nonumber\\
&&\quad+C2^{-2j/q}c_j(T)^\f1q\bigl(\|\dv \Delta_jG\|_{L^1_t(L^p)}+\|\Delta_jF_1\|_{L^1_t(L^p)}+
\|\dv[\Delta_j,\overline{\nu}]\na d\|_{L^1_t(L^p)}\bigr)\nonumber\\
&&\quad+C2^{-2j/q}c_j(T)^\f1q\bigl(\|\textrm{curl} \Delta_jG\|_{L^1_t(L^p)}+\|\Delta_jF_2\|_{L^1_t(L^p)}+\|\dv[\Delta_j,\overline{\mu}]\na w\|_{L^1_t(L^p)}\bigr),
\een
with $c_j(T)=1-e^{c_p2^{2j}T}$. Notice that
\beno
2^j\|\Delta_ju\|_p\sim \|\Delta_jd\|_p+\|\Delta_jw\|_p, \quad e_j(T)\le \om_j(T),
\eeno
which together with (\ref{eq:localcurl-dvest}) implies that
\ben\label{eq:momeninfty}
&&\|u\|_{\widetilde{L}^q_T(\dot B^{s-1+2/q}_{p,1})}
\le C\Bigl(\|u_0\|_{\dot B^{s-1}_{p,1}}+\|G\|_{\widetilde{L}^1_T(\dot B^{s-1}_{p,1})}+\|(F_1,F_2)\|_{\widetilde{L}^1_T(\dot B^{s-2}_{p,1})}\Bigr)\nonumber\\
&&\qquad+C\sum_j2^{j(s-2)}\bigl(\|\dv[\Delta_j,\overline{\nu}]\na
d\|_{L^1_T(L^p)}+\|\dv[\Delta_j,\overline{\mu}]\na
w\|_{L^1_T(L^p)}\bigr), \een and with $c=c_p$ in the definition of
$e_k(t)$, \ben\label{eq:momenone}
&&\|u\|_{\widetilde{L}^1_T(\dot B^{s-1}_{p,1})}+\|u\|_{\widetilde{L}^2_T(\dot B^{s}_{p,1})}\nonumber\\
&&\le C\bigl(\|u_0\|_{\dot B^{s-1}_{p,1}(\om)}+\|G\|_{\widetilde{L}^1_T(\dot B^{s-1}_{p,1}(\om))}+\|(F_1,F_2)\|_{\widetilde{L}^1_T(\dot B^{s-2}_{p,1}(\om))}\bigr)\nonumber\\
&&\quad+C\sum_j2^{j(s-2)}\om_j(T)\bigl(\|\dv[\Delta_j,\overline{\nu}]\na
d\|_{L^1_T(L^p)}+\|\dv[\Delta_j,\overline{\mu}]\na
w\|_{L^1_T(L^p)}\bigr). \een

First of all, we deal with the right hand side of (\ref{eq:momeninfty}). From Lemma \ref{Lem:binesti} and \ref{Lem:nonesti},
we infer that
\ben\label{eq:F1}
\|F_1\|_{\widetilde{L}^1_T(\dot B^{s-2}_{p,1})}&\le& C\bigl(\|\na \overline{\mu}\cdot\na u\|_{\widetilde{L}^1_T(\dot B^{s-1}_{p,1})}
+\|\na (\overline{\lambda}+\overline{\mu})d\|_{\widetilde{L}^1_T(\dot B^{s-1}_{p,1})}\bigr)\nonumber\\
&\le&
C\bigl(\|\overline{\mu}-\overline{\mu}(\underline{\rho})\|_{\widetilde{L}^\infty_T(\dot
B^{ \frac{N}{p}}_{p,1})}+
\|\overline{\lambda}-\overline{\lambda}(\underline{\rho})\|_{\widetilde{L}^\infty_T(\dot B^{ \frac{N}{p}}_{p,1})}\bigr)\|u\|_{\widetilde{L}^1_T(\dot B^{s+1}_{p,1})}\nonumber\\
&\le& CA(T)\|\rho-\underline{\rho}\|_{\widetilde{L}^\infty_T(\dot
B^{ \frac{N}{p}}_{p,1})}\|u\|_{\widetilde{L}^1_T(\dot
B^{s+1}_{p,1})}. \een Similarly, we have \ben\label{eq:F2}
\|F_2\|_{\widetilde{L}^1_T(\dot B^{s-2}_{p,1})} \le
CA(T)\|\rho-\underline{\rho}\|_{\widetilde{L}^\infty_T(\dot B^{
\frac{N}{p}}_{p,1})}\|u\|_{\widetilde{L}^1_T(\dot B^{s+1}_{p,1})}.
\een While, we write \beno [\Delta_j,\overline{\nu}]\na
d=[\Delta_j,\overline{\nu}-\overline{\nu}(\underline{\rho})]\na d
=\Delta_j((\overline{\nu}-\overline{\nu}(\underline{\rho}))\na
d)-(\overline{\nu}-\overline{\nu}(\underline{\rho}))\Delta_j\na d,
\eeno then by Lemma \ref{Lem:binesti}, Lemma \ref{Lem:nonesti} we
get for $p\in[1,N]$\beno
&&\sum_j2^{j(s-2)}\|\dv\Delta_j((\overline{\nu}-\overline{\nu}(\underline{\rho}))\na
d)\|_{L^1_T(L^p)} \le
CA(T)\|\rho-\underline{\rho}\|_{\widetilde{L}^\infty_T(\dot B^{
\frac{N}{p}}_{p,1})}\|u\|_{\widetilde{L}^1_T(\dot B^{s+1}_{p,1})},\\
&&\sum_j2^{j(s-2)}\Big(\|\overline{\nu}-\overline{\nu}(\underline{\rho})\|_{L^\infty_T(L^\infty)}\|\Delta_j\dv\na
d\|_{L^1_T(L^p)}+\|\dv(\overline{\nu}-\overline{\nu}(\underline{\rho}))\|_{L^\infty_T(L^N)}\|\Delta_j\na
d\|_{L^1_T(L^{\frac{pN}{N-p}})}\Big)\quad\\&&\quad\le
CA(T)\|\rho-\underline{\rho}\|_{\widetilde{L}^\infty_T(\dot B^{
\frac{N}{p}}_{p,1})}\|u\|_{\widetilde{L}^1_T(\dot B^{s+1}_{p,1})},
\eeno
which imply that
\ben\label{eq:right3}
\sum_{j}2^{j(s-2)}\|\dv[\Delta_j,\overline{\nu}]\na
d\|_{L^1_T(L^p)}\le
CA(T)\|\rho-\underline{\rho}\|_{\widetilde{L}^\infty_T(\dot B^{
\frac{N}{p}}_{p,1})}\|u\|_{\widetilde{L}^1_T(\dot
B^{s+1}_{p,1})}.\een Similarly, we have \ben\label{eq:right4}
\sum_{j}2^{j(s-2)}\|\dv[\Delta_j,\overline{\mu}]\na w\|_{L^1_T(L^p)}
\le CA(T)\|\rho-\underline{\rho}\|_{\widetilde{L}^\infty_T(\dot B^{
\frac{N}{p}}_{p,1})}\|u\|_{\widetilde{L}^1_T(\dot B^{s+1}_{p,1})}.
\een Then the first inequality of Proposition \ref{Prop:momenequ}
(a) can be deduced from (\ref{eq:momeninfty}) and
(\ref{eq:F1})-(\ref{eq:right4}). On the other hand, using Lemma
\ref{Lem:binesti} and \ref{Lem:nonesti}, we also have \beno
\|F_1\|_{\widetilde{L}^1_T(\dot
B^{s-2}_{p,1})}+\|F_2\|_{\widetilde{L}^1_T(\dot B^{s-2}_{p,1})} \le
CA(T)\|\rho-\underline{\rho}\|_{\widetilde{L}^\infty_T(\dot B^{
\frac{N}{p}+1}_{p,1})}\|u\|_{\widetilde{L}^1_T(\dot B^{s}_{p,1})},
\eeno and by Lemma \ref{Lem:commu}, \beno
&&\sum_j2^{j(s-2)}\bigl(\|\dv[\Delta_j,\overline{\nu}]\na d\|_{L^1_T(L^p)}+\|\dv[\Delta_j,\overline{\mu}]\na w\|_{L^1_T(L^p)}\bigr)\\
&& \quad\le
CA(T)\|\rho-\underline{\rho}\|_{\widetilde{L}^\infty_T(\dot B^{
\frac{N}{p}+1}_{p,1})}\|u\|_{\widetilde{L}^1_T(\dot B^{s}_{p,1})},
\eeno which together with (\ref{eq:momeninfty}) lead to the second
inequality of Proposition \ref{Prop:momenequ} (a).

Next, we deal with the right hand side of (\ref{eq:momenone}).
From Proposition \ref{Prop:binweight} and \ref{Prop:nonweight}, it follows that
\ben\label{eq:F1-w}
\|F_1\|_{\widetilde{L}^1_T(\dot B^{s-2}_{p,1}(\om))}&\le& C\bigl(\|\na \overline{\mu}\cdot\na u\|_{\widetilde{L}^1_T(\dot B^{s-1}_{p,1}(\om))}
+\|\na (\overline{\lambda}+\overline{\mu})d\|_{\widetilde{L}^1_T(\dot B^{s-1}_{p,1}(\om))}\bigr)\nonumber\\
&\le&
C\bigl(\|\overline{\mu}-\overline{\mu}(\underline{\rho})\|_{\widetilde{L}^\infty_T(\dot
B^{ \frac{N}{p}}_{p,1}(\om))}+
\|\overline{\lambda}-\overline{\lambda}(\underline{\rho})\|_{\widetilde{L}^\infty_T(\dot B^{ \frac{N}{p}}_{p,1}(\om))}\bigr)\|u\|_{\widetilde{L}^1_T(\dot B^{s+1}_{p,1})}\nonumber\\
&\le& CA(T)\|\rho-\underline{\rho}\|_{\widetilde{L}^\infty_T(\dot
B^{ \frac{N}{p}}_{p,1}(\om))}\|u\|_{\widetilde{L}^1_T(\dot
B^{s+1}_{p,1})}. \een Similarly, we have \ben\label{eq:F2-w}
\|F_2\|_{\widetilde{L}^1_T(\dot B^{s-2}_{p,1}(\om))} \le
CA(T)\|\rho-\underline{\rho}\|_{\widetilde{L}^\infty_T(\dot B^{
\frac{N}{p}}_{p,1}(\om))}\|u\|_{\widetilde{L}^1_T(\dot
B^{s+1}_{p,1})}. \een

Notice that
\beno
[\Delta_j,\overline{\nu}]\na d=[\Delta_j,\overline{\nu}-\overline{\nu}(\underline{\rho})]\na d,\quad
[\Delta_j,\overline{\mu}]\na w=[\Delta_j,\overline{\mu}-\overline{\mu}(\underline{\rho})]\na w,
\eeno
which together with Lemma \ref{Lem:timecommweight} and Proposition \ref{Prop:nonweight} ensures that
\ben\label{eq:comm-w}
&&\sum_j2^{j(s-2)}\omega_j(T)\bigl(\|\dv[\Delta_j,\overline{\nu}]\na d\|_{L^1_T(L^p)}+\|\dv[\Delta_j,\overline{\mu}]\na w\|_{L^1_T(L^p)}\bigr)\nonumber\\
&&\quad\le
C\bigl(\|\overline{\nu}-\overline{\nu}(\underline{\rho})\|_{\widetilde{L}^\infty_T(\dot
B^{ \frac{N}{p}}_{p,1}(\om))}
+\|\overline{\mu}-\overline{\mu}(\underline{\rho})\|_{\widetilde{L}^\infty_T(\dot
B^{ \frac{N}{p}}_{p,1}(\om))}\bigr)
\|u\|_{\widetilde{L}^1_T(\dot B^{s+1}_{p,1})}\nonumber\\
&&\quad\le
A(T)\|\rho-\underline{\rho}\|_{\widetilde{L}^\infty_T(\dot B^{
\frac{N}{p}}_{p,1}(\om))}\|u\|_{\widetilde{L}^1_T(\dot
B^{s+1}_{p,1})}. \een Summing up (\ref{eq:momenone}) and
(\ref{eq:F1-w})-(\ref{eq:comm-w}), we obtain the  inequality
of Proposition \ref{Prop:momenequ} (b).\ef

\begin{Lemma}\label{Lem:timecommweight}
Let $p\in [1,N]$ and $s\in (-N\min(\frac 1p,\frac1{p'}),\frac Np]$.
Assume that $f\in \widetilde{L}^\infty_T(\dot B^{
\frac{N}{p}}_{p,1}(\om))$ and $g\in \widetilde{L}^1_T(\dot
B^{s+1}_{p,1})$. Then  there holds \beno
\sum_j2^{j(s-1)}\omega_j(T)\|{\rm div}[\Delta_j,f]\na
g\|_{L^1_T(L^p)}\le C\|f\|_{\widetilde{L}^\infty_T(\dot B^{
\frac{N}{p}}_{p,1}(\om))}\|g\|_{\widetilde{L}^1_T(\dot
B^{s+1}_{p,1})}. \eeno
\end{Lemma}

\no{\bf Proof.}\, Using the Bony's decomposition, we write
\beno
[f,\Delta_j]\p_kg&=&[T_{f},\Delta_j]\p_kg+T_{\p_k\Delta_jg}f+R(f,\p_k\Delta_jg)\\
&&-\Delta_j(T_{\p_kg}f)-\Delta_jR(f,\p_kg). \eeno Using Lemma
\ref{Lem:bonyweighttimespace} (a) and (c) with $s_1= \frac{N}{p}$ and
$s_2=s$, we get \beno
&&\sum_{j}\om_j(T)2^{j(s-1)}\|\dv\Delta_j(T_{\p_kg}f)\|_{L^1_T(L^p)}\le C\|f\|_{\widetilde{L}^\infty_T(\dot B^{ \frac{N}{p}}_{p,1}(\om))}\|g\|_{\widetilde{L}^1_T(\dot B^{s+1}_{p,1})},\\
&&\sum_{j}\om_j(T)2^{j(s-1)}\|\dv\Delta_jR(f,\p_kg)\|_{L^1_T(L^p)}\le
C\|f\|_{\widetilde{L}^\infty_T(\dot B^{
\frac{N}{p}}_{p,1}(\om))}\|g\|_{\widetilde{L}^1_T(\dot
B^{s+1}_{p,1})}. \eeno Thanks to the proof of Lemma
\ref{Lem:commweight}, we have \beno T_{\p_k\Delta_jg}'f\triangleq
T_{\p_k\Delta_jg}f+R(f,\p_k\Delta_jg)=\sum_{j'\ge
j-2}S_{j'+2}\Delta_j\p_kg\Delta_{j'}f, \eeno then we get by Lemma
\ref{Lem:Bernstein} and (\ref{weightprop}) that \beno
&&\sum_{j}\om_j(T)2^{j(s-1)}\|\dv T_{\p_k\Delta_jg}'f\|_{L^1_T(L^p)}\\
&&\le C\sum_{j}\om_j(T)2^{js}\Big(\|\Delta_j\na g\|_{L^1_T(L^\infty)}\sum_{j'\ge j-2}\|\Delta_{j'}f\|_{L^1_T(L^p)}+
\|\Delta_j g\|_{L^1_T(L^\infty)}\sum_{j'\ge j-2}2^{j'}\|\Delta_{j'}f\|_{L^1_T(L^p)}\Big)\\
&&\le C\sum_{j}2^{j(s+ \frac{N}{p})}\|\Delta_jg\|_{L^1_T(L^p)}\sum_{j'\ge j-2}(2^j+2^{j'})\om_{j'}(T)\|\Delta_{j'}f\|_{L^1_T(L^p)}\\
&&\le C\|f\|_{\widetilde{L}^\infty_T(\dot B^{
\frac{N}{p}}_{p,1}(\om))}\|g\|_{\widetilde{L}^1_T(\dot
B^{s+1}_{p,1})}. \eeno Set $h(x)=(\cF^{-1}\varphi)(x)$. Thanks to
the proof of Lemma \ref{Lem:commweight}, we have \beno [T_{f},
\Delta_j]\pa_k g
&=&\sum_{|j'-j|\le4}2^{(N+1)j}\int_{\R^N}\int_0^1y\cdot\na
S_{j'-1}f(x-\tau
y)d\tau\pa_kh(2^jy)\Delta_{j'}g(x-y)dy\nonumber\\&&\qquad\quad+
2^{Nj}\int_{\R^N}h(2^j(x-y))\pa_kS_{j'-1}f(y)\Delta_{j'}g(y)dy ,
\eeno from which and a similar argument of Lemma
\ref{Lem:bonyweighttimespace} (b), we infer that \beno
\sum_j\om_j(T)2^{j(s-1)}\|\dv[T_{f}, \Delta_j]\pa_k
g\|_{L^1_T(L^p)}\le C\|f\|_{\widetilde{L}^\infty_T(\dot B^{
\frac{N}{p}}_{p,1}(\om))}\|g\|_{\widetilde{L}^1_T(\dot
B^{s+1}_{p,1})}. \eeno

Summing up all the above estimates, we conclude the proof of Lemma \ref{Lem:timecommweight}.\ef

To prove the uniqueness of the solution, we also need the following proposition.

\begin{Proposition}\label{Prop:momenequ-end} Let $p\in [2,N]$.
Assume that $G\in L^1_T(\dot B^{- \frac{N}{p}}_{p,\infty}), u_0\in
\dot B^{- \frac{N}{p}}_{p,\infty}$, and $\rho-\underline{\rho}\in
L^\infty_T(\dot B^{ \frac{N}{p}}_{p,1})$. Let $u$ be a solution of
(\ref{eq:linearmomen}). Then there holds \beno
&&\|u\|_{\widetilde{L}^1_T(\dot B^{- \frac{N}{p}+2}_{p,\infty})}+\|u\|_{\widetilde{L}^2_T(\dot B^{- \frac{N}{p}+1}_{p,\infty})}\\
&&\le C\Bigl(\|u_0\|_{\dot B^{-
\frac{N}{p}}_{p,\infty}}+\|G(\tau)\|_{\widetilde{L}^1_T(\dot B^{-
\frac{N}{p}}_{p,\infty}(\om))}
+A(T)\|\rho-\underline{\rho}\|_{\widetilde{L}^\infty_T(\dot B^{
\frac{N}{p}}_{p,1}(\om))}\|u\|_{\widetilde{L}^1_T(\dot B^{-
\frac{N}{p}+2}_{p,\infty})}\Bigr). \eeno
\end{Proposition}

\noindent{\bf Proof.}\, We closely follow the proof of Proposition \ref{Prop:momenequ}.
From (\ref{eq:localcurl-dvest}), we infer that
\ben\label{eq:momenone-end}
&&\|u\|_{\widetilde{L}^1_T(\dot B^{- \frac{N}{p}+2}_{p,\infty})}+\|u\|_{\widetilde{L}^2_T(\dot B^{- \frac{N}{p}+1}_{p,\infty})}\nonumber\\
&&\le C\Bigl(\|u_0\|_{\dot B^{- \frac{N}{p}}_{p,\infty}(\om)}+\|G\|_{\widetilde{L}^1_T(\dot B^{- \frac{N}{p}}_{p,\infty}(\om))}+\|(F_1,F_2)\|_{\widetilde{L}^1_T(\dot B^{- \frac{N}{p}-1}_{p,\infty}(\om))}\Bigr)\nonumber\\
&&\quad+C\sup_{j\in \Z}2^{j(
-\frac{N}{p}-1)}\omega_j(T)\bigl(\|\dv[\Delta_j,\overline{\nu}]\na
d\|_{L^1_T(L^p)}+\|\dv[\Delta_j,\overline{\mu}]\na
w\|_{L^1_T(L^p)}\bigr). \een We use Proposition
\ref{Prop:binweight-end} to get \beno
\|(F_1,F_2)\|_{\widetilde{L}^1_T(\dot B^{-
\frac{N}{p}-1}_{p,\infty}(\om))} \le
CA(T)\|\rho-\underline{\rho}\|_{\widetilde{L}^\infty_T(\dot B^{
\frac{N}{p}}_{p,1}(\om))} \|u\|_{\widetilde{L}^1_T(\dot B^{-
\frac{N}{p}+2}_{p,\infty})}. \eeno From Lemma
\ref{Lem:timecommweight-end}, the second term on the right hand side
of (\ref{eq:momenone-end}) is bounded by \beno
CA(T)\|\rho-\underline{\rho}\|_{\widetilde{L}^\infty_T(\dot B^{
\frac{N}{p}}_{p,1}(\om))} \|u\|_{\widetilde{L}^1_T(\dot B^{-
\frac{N}{p}+2}_{p,\infty})}. \eeno

This completes the proof of Proposition \ref{Prop:momenequ-end}.\ef

\begin{Lemma}\label{Lem:timecommweight-end}
Let $p\in [1,N]$. Assume that $f\in \widetilde{L}^\infty_T(\dot B^{
\frac{N}{p}}_{p,1}(\om))$ and $g\in \widetilde{L}^1_T(\dot B^{-
\frac{N}{p}+1}_{p,\infty})$. Then  there holds \beno \sup_{j\in \Z}2^{j(-
\frac{N}{p}-1)}\omega_j(T)\|{\rm div} [\Delta_j,f]\na
g\|_{L^1_T(L^p)}\le C\|f\|_{\widetilde{L}^\infty_T(\dot B^{
\frac{N}{p}}_{p,1}(\om))}\|g\|_{\widetilde{L}^1_T(\dot B^{-
\frac{N}{p}+1}_{p,\infty})}. \eeno
\end{Lemma}

The proof of Lemma \ref{Lem:timecommweight-end} is very similar to that of Lemma \ref{Lem:timecommweight}. Here we omit its proof.

\setcounter{equation}{0}
\section{The proof of existence}
We set \beno a(t,x)=\f {\rho(t,x)-\overline{\rho}_0}
{\overline{\rho}_0},\quad \overline{\mu}(\rho)=\f
{\mu(\rho)}{\rho},\quad
\overline{\lambda}(\rho)=\f
{\lambda(\rho)}{\rho}.
\eeno
Then the system (\ref{eq:cNS}) reads \ben\label{eq:linearcNS} \left\{
\begin{array}{ll}
\p_ta+u\cdot\na a=F,\\
\p_tu-\dv(\overline{\mu}\na u)
-\na((\overline{\lambda}+\overline{\mu})\dv\,u)=G, \\
(a,u)|_{t=0}=(a_0,u_0),
\end{array}
\right. \een with $a_0=\f {\rho_0(x)-\overline{\rho}_0}
{\overline{\rho}_0}$ and \beno
&&F(a,u)=-(1+a)\dv\, u,\\
&&G(a,u)=-u\cdot\na u+\f{\overline{\rho}_0P'(\rho)} {\rho}\na a+\f
{\mu(\rho)} {\rho^2}\na \rho\cdot\na u+\f
{\mu(\rho)+\lambda(\rho)}{\rho^2}\na \rho\dv \,u. \eeno

\textbf{Step 1}.\, The approximate solution sequence\vspace{0.1cm}

We smooth out the data as follows:
\beno
a^n_0=S_{n+N}a_0,\quad u_0^n=S_nu_0,
\eeno
where $N\in\Z$ is chosen such that
\ben\label{assu:N}
\overline{\rho}_0(1+a^n_0(x))\ge \f {3} 4c_0.
\een

A standard linearized argument (as in the  proof of Theorem 4.2 in
\cite{Dan-cpde}) will ensure that the system (\ref{eq:linearcNS})
with the smooth data $(a_0^n,u_0^n)$ has a solution $(a^n,u^n)$ on a
time interval $[0,T_n]$ for some $T_n>0$ such that
\begin{equation}\label{eq:appsolutionspace}
\begin{split}
&a^n\in C([0,T_n];\dot B^{ \frac{N}{p}}_{p,1}\cap \dot B^{ \frac{N}{p}+1}_{p,1}) \quad \textrm{and} \\
&u^n\in C([0,T_n];\dot B^{ \frac{N}{p}-1}_{p,1}\cap \dot B^{
\frac{N}{p}}_{p,1})\cap L^1([0,T_n];\dot B^{
\frac{N}{p}+1}_{p,1}\cap \dot B^{ \frac{N}{p}+2}_{p,1}).
 \end{split}
\end{equation}
In what follows, we also denote by $T_n$ the maximal lifespan of the solution $(a^n,u^n)$.\vspace{0.2cm}

\textbf{Step 2}.\, Uniform estimates\vspace{0.1cm}

Let $E_0:=\|a_0\|_{\dot B^{ \frac{N}{p}}_{p,1}}+\|u_0\|_{\dot B^{
\frac{N}{p}-1}_{p,1}}$ and $T\in (0,T_n)$. We assume that the
solutions $(a^n,u^n)$ satisfies the following inequalities for some
positive constants $c_1, C_0, A_0$ and $\eta$(to be determined
later):\vspace{0.1cm}

\no(\cH 1)\, $\overline{\rho}_0(1+a_0^n(t,x))\ge \f {c_0} 2$ for any
$(t,x)\in [0,T]\times \R^N$;\vspace{0.1cm}

\no(\cH 2)\, $\overline{\mu}^n(t,x)\ge c_1,\overline{\lambda}^n(t,x)+2\overline{\mu}^n(t,x)\ge c_1$ for any $(t,x)\in [0,T]\times \R^N$;\vspace{0.1cm}

\no(\cH 3)\, $\|a^n\|_{\widetilde{L}^\infty_T(\dot B^{
\frac{N}{p}}_{p,1})}+\|u^n\|_{\widetilde{L}^\infty_T(\dot B^{
\frac{N}{p}-1}_{p,1})}\le C_0E_0;$\vspace{0.1cm}

\no(\cH 4)\, $\|a^n\|_{\widetilde{L}^\infty_T(\dot B^{
\frac{N}{p}}_{p,1}(\om))}\le A_0\eta,
\,\|u^n\|_{\widetilde{L}^1_T(\dot B^{ \frac{N}{p}+1}_{p,1})}
+\|u^n\|_{\widetilde{L}^2_T(\dot B^{ \frac{N}{p}}_{p,1})}\le
\eta.$\vspace{0.1cm}

In what follows, we will show that if the conditions (\cH 1) to (\cH 4) are satisfied for some $T>0$,
then they are actually satisfied with strict inequalities. Since all those conditions depend continuously
on the time variable and are satisfied initially, a standard bootstrap argument will ensure that
(\cH 1) to (\cH 4) are indeed satisfied for $T$.

First of all, we get by Proposition \ref{Prop:transport} that
\ben\label{eq:dens1} \|a^n\|_{\widetilde{L}^\infty_T(\dot{B}^{
\frac{N}{p}}_{p,1})}\le e^{CV^n(T)}\bigl( \|a_0\|_{\dot{B}^{
\frac{N}{p}}_{p,1}}+\|F^n\|_{\widetilde{L}^1_T(\dot{B}^{
\frac{N}{p}}_{p,1})}\bigr), \een and by Proposition
\ref{Prop:momenequ}, we have \ben\label{eq:veloc1}
\|u^n\|_{\widetilde{L}^\infty_T(\dot B^{\frac{N}{p}-1}_{p,1})}
&\le& C\|u_0\|_{\dot B^{ \frac{N}{p}-1}_{p,1}}+C\|G^n\|_{\widetilde{L}^1_T(\dot B^{ \frac{N}{p}-1}_{p,1})}\nonumber\\
&&+CA^n(T)\|a^n\|_{\widetilde{L}^\infty_T(\dot B^{
\frac{N}{p}}_{p,1})}\|u^n\|_{\widetilde{L}^1_T(\dot B^{
\frac{N}{p}+1}_{p,1})}, \een where $V^n(t)=\int_0^t\|\nabla
u^n(\tau)\|_{\dot{B}^{\frac{N}{p}}_{p,1}}d\tau$ and
$A^n(T)=\bigl(1+\|a^n\|_{\widetilde{L}^\infty_T(\dot B^{
\frac{N}{p}}_{p,1})}\bigr)^{[ \frac{N}{p}]+3}$. For $F^n$, we apply
Lemma \ref{Lem:binesti} to get \beno
\|F^n\|_{\widetilde{L}^1_T(\dot{B}^{ \frac{N}{p}}_{p,1})}\le
\|u^n\|_{\widetilde{L}^1_T(\dot{B}^{ \frac{N}{p}+1}_{p,1})}+
C\|a^n\|_{\widetilde{L}^\infty_T(\dot{B}^{
\frac{N}{p}}_{p,1})}\|u^n\|_{\widetilde{L}^1_T(\dot{B}^{
\frac{N}{p}+1}_{p,1})}, \eeno and for $G^n$, we use Lemma
\ref{Lem:binesti} and \ref{Lem:nonesti} to get \beno
\|G^n\|_{\widetilde{L}^1_T(\dot B^{ \frac{N}{p}-1}_{p,1})}
&\le& C\|u^n\|_{\widetilde{L}^\infty_T(\dot{B}^{ \frac{N}{p}-1}_{p,1})}\|u^n\|_{\widetilde{L}^1_T(\dot{B}^{ \frac{N}{p}+1}_{p,1})}\\
&&+CA^n(T)\|a^n\|_{\widetilde{L}^\infty_T(\dot{B}^{
\frac{N}{p}}_{p,1})} \bigl(T+\|u^n\|_{\widetilde{L}^1_T(\dot{B}^{
\frac{N}{p}+1}_{p,1})}\bigr). \eeno Plugging the above two estimates
into (\ref{eq:dens1}) and (\ref{eq:veloc1}), we obtain
\ben\label{eq:denvel-infty}
\|a^n\|_{\widetilde{L}^\infty_T(\dot{B}^{
\frac{N}{p}}_{p,1})}+\|u^n\|_{\widetilde{L}^\infty_T(\dot
B^{\frac{N}{p}-1}_{p,1})} \le
C_1e^{C_1V^n(T)}(E_0+(C_0E_0+1)\eta)+C_1A^n(T)C_0E_0(T+\eta). \een

Next, we get by Proposition \ref{Prop:transportweight} that
\ben\label{eq:dens2} \|a^n\|_{\widetilde{L}^\infty_T(\dot{B}^{
\frac{N}{p}}_{p,1}(\om))}\le e^{CV^n(T)}\bigl( \|a_0\|_{\dot{B}^{
\frac{N}{p}}_{p,1}(\om)}+\|F^n\|_{\widetilde{L}^1_T\dot{B}^{
\frac{N}{p}}_{p,1}(\om)}\bigr), \een and by Proposition
\ref{Prop:momenequ}, we have \ben\label{eq:veloc2}
&&\|u^n\|_{\widetilde{L}^1_T(\dot
B^{\frac{N}{p}+1}_{p,1})}+\|u^n\|_{\widetilde{L}^2_T(\dot
B^{\frac{N}{p}}_{p,1})}\nonumber\\&& \le C\Bigl(\|u_0\|_{\dot
B^{\frac{N}{p}-1}_{p,1}(\om))} +\|G^n\|_{\widetilde{L}^1_T(\dot
B^{\frac{N}{p}-1}_{p,1}(\om))}
+A^n(T)\|a^n\|_{\widetilde{L}^\infty_T(\dot B^{
\frac{N}{p}}_{p,1}(\om))}\|u^n\|_{\widetilde{L}^1_T(\dot
B^{\frac{N}{p}+1}_{p,1})}\Bigr). \een For $F^n$, we have \beno
\|F^n\|_{\widetilde{L}^1_T(\dot{B}^{ \frac{N}{p}}_{p,1}(\om))}\le
2\|F^n\|_{\widetilde{L}^1_T(\dot{B}^{ \frac{N}{p}}_{p,1})}\le C(1+
\|a^n\|_{\widetilde{L}^\infty_T(\dot{B}^{
\frac{N}{p}}_{p,1})})\|u^n\|_{\widetilde{L}^1_T(\dot{B}^{
\frac{N}{p}+1}_{p,1})}, \eeno and for $G^n$, we use Proposition
\ref{Prop:binweight} and \ref{Prop:nonweight} to get \beno
\|G^n\|_{\widetilde{L}^1_T(\dot B^{ \frac{N}{p}-1}_{p,1}(\om))} \le
C\|u^n\|_{\widetilde{L}^2_T(\dot{B}^{ \frac{N}{p}}_{p,1})}^2
+CA^n(T)\|a^n\|_{\widetilde{L}^\infty_T(\dot{B}^{
\frac{N}{p}}_{p,1}(\om))}
\bigl(T+\|u^n\|_{\widetilde{L}^1_T(\dot{B}^{
\frac{N}{p}+1}_{p,1})}\bigr). \eeno Plugging the above two estimates
into (\ref{eq:dens2}) and (\ref{eq:veloc2}), we obtain
\ben\label{eq:dens-weight}
&&\|a^n\|_{\widetilde{L}^\infty_T(\dot{B}^{
\frac{N}{p}}_{p,1}(\om))}\le e^{C_1V^n(T)}\bigl(
\|a_0\|_{\dot{B}^{ \frac{N}{p}}_{p,1}(\om)}+C_2(1+C_0E_0)\eta\bigr),\\
&&\|u^n\|_{\widetilde{L}^1_T(\dot B^{\frac
Np+1}_{p,1})}+\|u^n\|_{\widetilde{L}^2_T(\dot B^{\frac Np}_{p,1})}
\le C_3\Bigl(\|u_0\|_{\dot B^{\frac Np-1}_{p,1}(\om))}
+\eta^2+A_0A^n(T)\eta(T+\eta)\Bigr).\quad\label{eq:veloc-onetwo}
\een

According to the definition of $V^n$ and $A^n$, we have \beno
V^n(T)\le \eta,\quad A^n(T)\le (1+C_0E_0)^{[ \frac{N}{p}]+3}. \eeno

Let $C_0=4C_1$ and $A_0=2C_2(1+C_0E_0)$. Then we take $\eta$ small
enough such that \begin{equation} \label{assu:eta}
\begin{split}
&e^{C_1\eta}<\frac 3 2,\,(C_0E_0+1)\eta \le E_0,\, C_1(C_0E_0+1)^{[
\frac{N}{p}]+3}\eta \le \f{1} {16}, \\ & C_3\eta\le \f16, \,
C_3A_0(1+C_0E_0)^{[ \frac{N}{p}]+3}\eta\le \f16.
\end{split}
\end{equation}
Next, we take $T$ small enough such that \ben\label{assu:time1}
C_1(C_0E_0+1)^{[ \frac{N}{p}]+3}T \le \f{1}{16}, \,
C_3A_0(1+C_0E_0)^{[ \frac{N}{p}]+3}T\le \f16. \een and note that
$\om_k(0)=0$ and $(a_0,u_0)\in \dot{B}^{ \frac{N}{p}}_{p,1}\times
\dot{B}^{ \frac{N}{p}-1}_{p,1}$, we can also take $T$ small enough
such that \ben\label{assu:time2} \|a_0\|_{\dot{B}^{
\frac{N}{p}}_{p,1}(\om)}\le \f {A_0} {12}\eta,\quad
\|u_0\|_{\dot{B}^{ \frac{N}{p}-1}_{p,1}(\om)}\le \f \eta {6C_3}.
\een Then it follows from (\ref{eq:denvel-infty}) that \beno
\|a^n\|_{\widetilde{L}^\infty_T(\dot{B}^{
\frac{N}{p}}_{p,1})}+\|u^n\|_{\widetilde{L}^\infty_T(\dot B^{\frac
Np-1}_{p,1})} \le \f78C_0E_0, \eeno and from (\ref{eq:dens-weight})
and (\ref{eq:veloc-onetwo}), we infer that \beno
\|a^n\|_{\widetilde{L}^\infty_T(\dot B^{
\frac{N}{p}}_{p,1}(\om))}\le \f78A_0\eta,
\,\|u^n\|_{\widetilde{L}^1_T(\dot B^{ \frac{N}{p}+1}_{p,1})}
+\|u^n\|_{\widetilde{L}^2_T(\dot B^{ \frac{N}{p}}_{p,1})}\le
\f23\eta, \eeno which ensure that (\cH 3) and (\cH 4) are satisfied
with strict inequalities for $T$ and $\eta$ satisfying
(\ref{assu:eta})-(\ref{assu:time2}).

Let $X_n(t,x)$ be a solution of \beno \f d {dt}
X_n(t,x)=u^n(t,X_n(t,x)),\quad X_n(0,x)=x, \eeno and we denote by
$X^{-1}_n(t,x)$ the inverse of $X_n(t,x)$. Then $a^n(t,x)$ can be
solved as \beno
a^n(t,x)=a_0^n(X^{-1}_n(t,x))+\int_0^tF^n(\tau,X_n(\tau,X^{-1}_n(t,x)))d\tau,
\eeno thus, we have \ben\label{eq:rho}
\overline{\rho}_0(1+a^n(t,x))=\rho_0^n(X^{-1}_n(t,x))+\bar{\rho}_0\int_0^tF^n(\tau,X_n(\tau,X^{-1}_n(t,x)))d\tau.
\een On the other hand, we have \beno \|F^n\|_{L^1_T(L^\infty)}\le
\|\na u\|_{L^1_T(L^\infty)}(1+\|a\|_{L^\infty_T(L^\infty)}) \le
C_4(1+C_0E_0)\eta. \eeno We take $\eta$ such that \beno
C_4(1+C_0E_0)\eta<\f 1 8c_0. \eeno Then from (\ref{eq:rho}) and
(\ref{assu:N}), it follows that \beno
\overline{\rho}_0(1+a^n(t,x))\ge \f34c_0-\f18c_0\ge \f58c_0, \eeno
that is, (\cH1) is  satisfied with the strict inequality. Finally,
take \beno c_1=\f12\min(\inf_{
|\rho|\le\overline{\rho}_0(1+C_0E_0)}
\overline{\mu}(\rho),\inf_{
|\rho|\le\overline{\rho}_0(1+C_0E_0)}
(\overline{\lambda}(\rho)+2\overline{\mu}(\rho))),
\eeno which ensures that (\cH 2) is satisfied with strict
inequality.

Let $T^*$ be the supremum of all time $T$ such that (\ref{assu:time1}) and (\ref{assu:time2}) are satisfied. We need to
prove that $T_n\ge T^*$. If $T_n<T^*$, then we can prove that
\begin{equation}\label{eq:appsolutionspace1}
\begin{split}
&a^n\in \widetilde{L}^\infty(0,T_n;\dot B^{ \frac{N}{p}}_{p,1}\cap \dot B^{ \frac{N}{p}+1}_{p,1}) \quad \textrm{and} \\
&u^n\in \widetilde{L}^\infty(0,T_n;\dot B^{ \frac{N}{p}-1}_{p,1}\cap
\dot B^{ \frac{N}{p}}_{p,1})\cap L^1([0,T_n];\dot B^{
\frac{N}{p}+1}_{p,1}\cap \dot B^{ \frac{N}{p}+2}_{p,1}),
 \end{split}
\end{equation}
thus, the solution $(a^n,u^n)$ can be continued beyond $T^*$.
Indeed, from Proposition \ref{Prop:transport}, we have
\ben\label{eq:dens3} \|a^n\|_{\widetilde{L}^\infty_t(\dot{B}^{
\frac{N}{p}+1}_{p,1})}\le e^{CV^n(t)}\bigl( \|a_0^n\|_{\dot{B}^{
\frac{N}{p}+1}_{p,1}}+\int_0^t\|F^n(\tau)\|_{\dot{B}^{
\frac{N}{p}+1}_{p,1}}d\tau\bigr), \een and by Proposition
\ref{Prop:momenequ} (a), we have \ben\label{eq:veloc3}
&&\|u^n\|_{\widetilde{L}^\infty_T(\dot B^{ \frac{N}{p}}_{p,1})}+\|u^n\|_{\widetilde{L}^1_T(\dot B^{ \frac{N}{p}+2}_{p,1})}\nonumber\\
&&\le C\|u_0^n\|_{\dot B^{
\frac{N}{p}}_{p,1}}+C\|G^n\|_{\widetilde{L}^1_T(\dot B^{
\frac{N}{p}}_{p,1})} +CA^n(T)\|a^n\|_{\widetilde{L}^\infty_T(\dot
B^{ \frac{N}{p}+1}_{p,1})}\|u^n\|_{\widetilde{L}^1_T(\dot B^{
\frac{N}{p}+1}_{p,1})}, \een On the other hand, we use Lemma
\ref{Lem:binesti-1} and the embedding $\dot{B}^{ \frac{N}{p}}_{p,1}\hookrightarrow L^\infty$ to get \beno
\|F^n\|_{\dot{B}^{ \frac{N}{p}+1}_{p,1}}\le \|u^n\|_{\dot{B}^{
\frac{N}{p}+2}_{p,1}}+ C\|a^n\|_{\dot{B}^{ \frac{N}{p}}_{p,1}}
\|u^n\|_{\dot{B}^{ \frac{N}{p}+2}_{p,1}}+ C\|a^n\|_{\dot{B}^{
\frac{N}{p}+1}_{p,1}} \|u^n\|_{\dot{B}^{ \frac{N}{p}+1}_{p,1}},
\eeno and by Lemma \ref{Lem:binesti} and \ref{Lem:nonesti}, we have
\beno \|G^n\|_{\widetilde{L}^1_T(\dot B^{ \frac{N}{p}}_{p,1})}
&\le& C\|u^n\|_{\widetilde{L}^\infty_T\dot{B}^{ \frac{N}{p}}_{p,1}}\|u^n\|_{\widetilde{L}^1_T\dot{B}^{ \frac{N}{p}+1}_{p,1}}\\
&&+CA^n(T)\|a^n\|_{\widetilde{L}^\infty_T\dot{B}^{
\frac{N}{p}+1}_{p,1}} \bigl(T+\|u^n\|_{\widetilde{L}^1_T\dot{B}^{
\frac{N}{p}+1}_{p,1}}\bigr), \eeno which together with
(\ref{eq:dens3}), (\ref{eq:veloc3}) and (\cH3-\cH4) implies
(\ref{eq:appsolutionspace1}).

\vspace{0.2cm}

\textbf{Step 3}.\, Existence of a solution\vspace{0.1cm}

We will use a compact argument to prove that the approximate sequence $\{a^n,u^n\}_{n\in \N}$ tends to some function
$(a,u)$ which satisfies the system (\ref{eq:linearcNS}) in the sense of distribution.

Since $\{u^n\}$ is uniformly bounded in $L^1_T(\dot B^{
\frac{N}{p}+1}_{p,1})\cap L^\infty_T(\dot B^{
\frac{N}{p}-1}_{p,1})$, we get by the interpolation that
$\{u^n\}_{n\in \N}$ is also uniformly bounded in $L^q_T(\dot B^{
\frac{N}{p}-1+2/q}_{p,1})$ for any $q\in [1,\infty]$. By Lemma
\ref{Lem:binesti}, we have \beno
&&\|a^n\dv u^n\|_{\dot B^{ \frac{N}{p}-1}_{p,1}}\le C\|a^n\|_{\dot B^{ \frac{N}{p}}_{p,1}}\|u^n\|_{\dot B^{ \frac{N}{p}}_{p,1}},\\
&&\|u^n\cdot\na a^n\|_{\dot B^{ \frac{N}{p}-1}_{p,1}}\le
C\|a^n\|_{\dot B^{ \frac{N}{p}}_{p,1}}\|u^n\|_{\dot B^{
\frac{N}{p}}_{p,1}}, \eeno from which and the first equation of the
system (\ref{eq:linearcNS}), we infer that $\{\p_t a^n\}_{n\in \N}$
is uniformly bounded in $L^2_T(\dot B^{ \frac{N}{p}-1}_{p,1})$. On
the other hand, by Lemma \ref{Lem:binesti} and Lemma
\ref{Lem:nonesti}, we have \beno
&&\|u^n\cdot\na u^n\|_{\dot B^{ \frac{N}{p}-3/2}_{p,1}}\le C\|u^n\|_{\dot B^{ \frac{N}{p}-1}_{p,1}}\|u^n\|_{\dot B^{ \frac{N}{p}+1/2}_{p,1}},\\
&&\|\f{\overline{\rho}_0P'(\rho^n)} {\rho^n}\na a^n\|_{\dot B^{ \frac{N}{p}-1}_{p,1}}\le C(\|a^n\|_\infty)\|a^n\|_{\dot B^{ \frac{N}{p}}_{p,1}}
(1+\|a^n\|_{\dot B^{ \frac{N}{p}}_{p,1}}),\\
&&\|\dv(\overline{\mu}^n\na u^n)\|_{\dot B^{
\frac{N}{p}-3/2}_{p,1}}+\|\na((\overline{\lambda}^n+\overline{\mu}^n)\dv
u^n)\|_{\dot B^{ \frac{N}{p}-3/2}_{p,1}}
+\|G_1^n\|_{\dot B^{ \frac{N}{p}-3/2}_{p,1}}\\
&&\quad\le C(\|a^n\|_\infty)\|a^n\|_{\dot B^{
\frac{N}{p}}_{p,1}}(1+\|a^n\|_{\dot B^{ \frac{N}{p}}_{p,1}})\|u^n\|_{\dot B^{ \frac{N}{p}+1/2}_{p,1}}, \eeno
where $G_1^n\eqdefa \f{\mu(\rho^n)} {(\rho^n)^2}\na \rho^n\cdot\na
u^n+\f{\mu(\rho^n)+\lambda(\rho^n)}{(\rho^n)^2}\na \rho^n\dv \,u^n$.
Then, from the second equation of the system (\ref{eq:linearcNS}),
we infer that $\{\p_t u^n\}_{n\in \N}$ is uniformly bounded in
$L^\f43_T(\dot B^{ \frac{N}{p}-\f32}_{p,1}+\dot B^{
\frac{N}{p}-1}_{p,1})$.

Let $\{\chi_j\}_{j\in \N}$ be a sequence of smooth functions
supported in the ball $B(0,j+1)$ and equal to 1 on $B(0,j)$. The
above proof ensures that for any $j\in \N$, $\{\chi_j a^n\}_{n\in
\N}$ is uniformly bounded in $C^\f12([0,T];\dot B^{
\frac{N}{p}-1}_{p,1})$, and $\{\chi_ju^n\}_{n\in \N}$ is uniformly
bounded in $C^\f14([0,T];\dot B^{ \frac{N}{p}-\f32}_{p,1}+\dot B^{
\frac{N}{p}-1}_{p,1})$. Since the embedding  $\dot B^{
\frac{N}{p}-1}_{p,1}\cap \dot B^{ \frac{N}{p}}_{p,1}\hookrightarrow
\dot B^{ \frac{N}{p}-1}_{p,1}$ and $\dot B^{
\frac{N}{p}-3/2}_{p,1}\cap \dot B^{
\frac{N}{p}-1}_{p,1}\hookrightarrow \dot B^{ \frac{N}{p}-3/2}_{p,1}$
are locally compact, by applying Ascoli's theorem and Cantor's
diagonal process , there exists some function $(a,u)$ such that for
any $j\in\N$, \ben\label{eq:conver1}
\begin{split}
&\chi_j a^n\longrightarrow \chi_j a \quad \textrm{in}\quad C([0,T];\dot B^{ \frac{N}{p}-1}_{p,1}),\\
&\chi_j u^n\longrightarrow \chi_j u \quad \textrm{in}\quad
C([0,T];\dot B^{ \frac{N}{p}-\f32}_{p,1}),
\end{split}
\een
as $n$ tends to $\infty$(up to a subsequence). By the interpolation, we also have
\ben\label{eq:conver2}
\begin{split}
&\chi_j a^n\longrightarrow \chi_j a \quad \textrm{in}\quad C([0,T];\dot B^{ \frac{N}{p}-s}_{p,1}),\quad \forall\,0<s\le 1,\\
&\chi_j u^n\longrightarrow \chi_j u \quad \textrm{in}\quad
L^1([0,T];\dot B^{ \frac{N}{p}+s}_{p,1}),\quad \forall\,-\f32\le
s<1.
\end{split}
\een
Furthermore, we actually have
\ben\label{eq:bound}
(a,u)\in \widetilde{L}^\infty_T(\dot B^{\f {N} p}_{p,1})\otimes \Bigl(\widetilde{L}^\infty_T(\dot B^{\f {N}
p-1}_{p,1}) \cap L^1(0,T; \dot B^{\f {N} p+1}_{p,1})\Bigr),\quad \overline{\rho}_0(1+a(t,x))\ge \f {c_0} 2.
\een

With (\ref{eq:conver1})-(\ref{eq:bound}), it is a routine process to
verify that $(a,u)$ satisfies the system (\ref{eq:linearcNS}) in the
sense of distribution(see also \cite{Dan-inve}). Finally, following
the argument in \cite{Dan-inve}, we can show that $(a,u)\in
C([0,T];\dot B^{ \frac{N}{p}}_{p,1})\otimes C([0,T];\dot B^{
\frac{N}{p}-1}_{p,1})$.

\setcounter{equation}{0}
\section{The proof of uniqueness}

In this section, we prove the uniqueness of the solution. Assume that $(a^1,u^1)\in E^p_T$ and $(a^2,u^2)\in E^p_T$
are  two solutions of the system (\ref{eq:linearcNS}) with the same initial data. Without loss of generality, we may assume that
 $a^1$ satisfies
\beno \rho^1(t,x)=\overline{\rho}_0(1+a^1(t,x))\ge \f {c_0} 2. \eeno
for any $(t,x)\in [0,T]\times \R^N$. Since $a^2\in C([0,T];\dot B^{
\frac{N}{p}}_{p,1})$ and $\rho^2(0,x)\ge c_0$, there exists a
positive time $\widetilde{T}\in (0,T]$ such that \beno
\rho^2(t,x)=\overline{\rho}_0(1+a^2(t,x))\ge \f {c_0} 2. \eeno for
any $(t,x)\in [0,\widetilde{T}]\times \R^N$. Set $\delta a=a^1-a^2$
and $\delta u=u^1-u^2$. Then $(\delta a,\delta u)$ satisfies
\ben\label{eq:linearcNS-diffe} \left\{
\begin{array}{ll}
\p_t\delta a+u^2\cdot\na \delta a=\delta F-\delta u\cdot\na a^1,\\
\p_t\delta u-\dv(\overline{\mu}^1\na \delta u)
-\na((\overline{\lambda}^1+\overline{\mu}^1)\dv\,\delta u)=\delta G+\delta H, \\
(\delta a,\delta u)|_{t=0}=(0,0),
\end{array}
\right. \een
where
\beno
&&\delta F=F(a^1,u^1)-F(a^2,u^2),\quad \delta G=G(a^1,u^1)-G(a^2,u^2),\\
&&\delta H=\dv\bigl((\overline{\mu}^1-\overline{\mu}^2)\na u^2\bigr)
+\na\bigl((\overline{\lambda}^1-\overline{\lambda}^2+\overline{\mu}^1-\overline{\mu}^2)\dv\,u^2\bigr),
\eeno
with $\overline{\lambda}^i=\overline{\lambda}(a^i), \overline{\mu}^i=\overline{\mu}(a^i)$ for $i=1,2$.

In what follows, we set $U^i(t)=\int_0^t\|u^i(\tau)\|_{\dot{B}^{
\frac{N}{p}+1}_{p,1}}d\tau$ for $i=1,2$, and denote by $A_T$ a
constant depending  on $\|a^1\|_{\widetilde{L}^\infty_T(\dot B^{
\frac{N}{p}}_{p,1})}$ and $\|a^2\|_{\widetilde{L}^\infty_T(\dot B^{
\frac{N}{p}}_{p,1})}$. Due to the inclusion relation $E_T^p\subseteq
E^N_T$, it suffices to prove the uniqueness of the solution in
$E^N_T$. So, we take $p=N$ in the sequel.

We apply Proposition \ref{Prop:transport} to get for any $t\in [0,T]$,
\ben\label{eq:dens-diffe2}
\|\delta a(t)\|_{\dot{B}^{0}_{p,\infty}}\le e^{CU^2(t)}\int_0^t\bigl(\|\delta F(\tau)\|_{\dot{B}^{0}_{p,\infty}}+
\|\delta u\cdot\na a^1(\tau)\|_{\dot{B}^{0}_{p,\infty}}\bigr)d\tau.
\een
By Lemma \ref{Lem:binesti-end}, we have
\beno
&&\|\delta F(\tau)\|_{\dot{B}^{0}_{p,\infty}}+
\|\delta u\cdot\na a^1(\tau)\|_{\dot{B}^{0}_{p,\infty}}\\
&&\quad\le C\|u^2\|_{\dot{B}^{2}_{p,1}}\|\delta a\|_{\dot{B}^{0}_{p,\infty}}+C(1+\|a^1\|_{\dot{B}^{1}_{p,1}})\|\delta u\|_{\dot{B}^{1}_{p,1}}.
\eeno
Plugging it into (\ref{eq:dens-diffe2}), we get by Gronwall's inequality that
\ben\label{eq:dens-diffe3}
\|\delta a(t)\|_{\dot{B}^{0}_{p,\infty}}\le e^{CU^2(t)}\int_0^t(1+\|a^1\|_{\dot{B}^{1}_{p,1}})\|\delta u\|_{\dot{B}^{1}_{p,1}}d\tau.
\een
We use Proposition \ref{Prop:momenequ-end} to get for any $t\in [0,T]$,
\ben\label{eq:velo-diff2}
\|\delta u(t)\|_{\widetilde{L}^1_t(\dot B^{1}_{p,\infty})}+\|\delta u(t)\|_{\widetilde{L}^2_t(\dot B^{0}_{p,\infty})}
&\le& C\int_0^t
\bigl(\|\delta G(\tau)\|_{\dot B^{-1}_{p,\infty}(\om)}+\|\delta H(\tau)\|_{\dot B^{-1}_{p,\infty}}\bigr)d\tau\nonumber\\
&&+CA_T\|a^1\|_{\widetilde{L}^\infty_t(\dot B^{1}_{p,1}(\om))}\|\delta u\|_{\widetilde{L}^1_t(\dot B^{1}_{p,\infty})}.
\een
From Lemma \ref{Lem:binesti-end}, Proposition \ref{Prop:binweight-end} and \ref{Prop:nonweight},
we infer that for any $t\in [0,\widetilde{T}]$,
\ben
\|\delta H\|_{\dot B^{-1}_{p,\infty}}
&\le& A_T\|u^2\|_{\dot B^{2}_{p,1}}
\|\delta a\|_{\dot B^{0}_{p,\infty}},\label{eq:H-diff1}\\
\|\delta G\|_{\dot B^{-1}_{p,\infty}(\om)}
&\le& C\|(u^1,u^2)\|_{\dot B^{1}_{p,1}}\|\delta u\|_{\dot B^{0}_{p,\infty}}
+A_T\|a^1\|_{\dot B^{1}_{p,1}(\om)}
\|\delta u\|_{\dot B^{1}_{p,\infty}}\nonumber\\
&&+A_T(1+\|u^2\|_{\dot B^{2}_{p,1}})
\|\delta a\|_{\dot B^{0}_{p,\infty}}.\label{eq:G-diff1}
\een
We take $\widetilde{T}$ small enough such that
\beno
\|(u^1,u^2)\|_{\widetilde{L}^1_t(\dot B^{2}_{p,1})\cap \widetilde{L}^2_t(\dot B^{1}_{p,1})}
+\|(a^1,a^2)\|_{\widetilde{L}^\infty_t(\dot B^{1}_{p,1}(\om))}\ll 1.
\eeno
Thus, plugging (\ref{eq:H-diff1}) and (\ref{eq:G-diff1}) into (\ref{eq:velo-diff2}),
we infer that for any $t\in [0,\widetilde{T}]$,
\ben\label{eq:velo-diff3}
\|\delta u\|_{\widetilde{L}^1_t(\dot B^{1}_{p,\infty})}
\le A_T\int_0^t\bigl(1+\|(u^1,u^2)\|_{\dot B^{2}_{p,1}}\bigr)\|\delta a\|_{\dot B^{0}_{p,\infty}}d\tau.
\een
\begin{Lemma}\cite{Dan-PRSE}\label{Lem:loginequ} Let $s\in \R$. Then for any $1\le p,\rho\le+\infty$ and $0<\epsilon\le 1$,
we have
\beno
\|f\|_{\widetilde{L}^\rho_T(\dot{B}^s_{p,1})}\le
C\frac{\|f\|_{\widetilde{L}^\rho_T(\dot{B}^s_{p,\infty})}}{\epsilon}
\log\Bigl(e+\frac{\|f\|_{\widetilde{L}^\rho_T(\dot{B}^{s-\epsilon}_{p,\infty})}
+\|f\|_{\widetilde{L}^\rho_T(\dot{B}^{s+\epsilon}_{p,\infty})}}
{\|f\|_{\widetilde{L}^\rho_T(\dot{B}^s_{p,\infty})}}\Bigr).
\eeno
\end{Lemma}
From Lemma \ref{Lem:loginequ}, it follows that
\beno
\|\delta u\|_{L^1_t(\dot B^{1}_{p,1})}\le C\|\delta u\|_{\widetilde{L}^1_t(\dot{B}^1_{p,\infty})}
\log\Bigl(e+\frac{\|\delta u\|_{\widetilde{L}^1_t(\dot{B}^{0}_{p,\infty})}
+\|\delta u\|_{\widetilde{L}^1_t(\dot{B}^{2}_{p,\infty})}}
{\|\delta u\|_{\widetilde{L}^1_t(\dot{B}^1_{p,\infty})}}\Bigr),
\eeno
which together with (\ref{eq:dens-diffe3}) and (\ref{eq:velo-diff3}) yields that for any $t\in [0,\widetilde{T}]$,
\beno
\|\delta u\|_{\widetilde{L}^1_t(\dot B^{1}_{p,\infty})}
\le A_T\int_0^t\bigl(1+\|(u^1,u^2)\|_{\dot B^{2}_{p,1}}\bigr)
\|\delta u\|_{\widetilde{L}^1_\tau(\dot{B}^1_{p,\infty})}
\log\bigl(e+C_T\|\delta u\|_{\widetilde{L}^1_\tau(\dot B^{1}_{p,\infty})}^{-1}\bigr)d\tau,
\eeno
where $C_T=\|\delta u\|_{\widetilde{L}^1_T(\dot{B}^{0}_{p,\infty})}
+\|\delta u\|_{\widetilde{L}^1_T(\dot{B}^{2}_{p,\infty})}$.
Notice that $1+\|(u^1,u^2)(t)\|_{\dot B^{2}_{p,1}}$ is integrable on $[0,T]$, and
\beno
\int_0^1\f {dr} {r\log(e+C_Tr^{-1})}dr=+\infty,
\eeno
Osgood lemma applied concludes that $(\delta a,\delta u)=0$ on $[0,\widetilde{T}]$, and
a continuity argument ensures that $(a^1,u^1)=(a^2,u^2)$ on $[0,T]$.

\bigskip

\section*{Acknowledgements}
 Q. Chen and C. Miao were partially supported by the NSF of China under grant
No.10701012, No.10725102. Z. Zhang was partially supported by the
NSF of China under grant No.10601002.

\end{document}